\documentclass{ieot}

 \newtheorem{thm}{Theorem}[section]
 \newtheorem{cor}[thm]{Corollary}
 \newtheorem{lem}[thm]{Lemma}
 \newtheorem{prop}[thm]{Proposition}
 \theoremstyle{definition}
 
 \theoremstyle{remark}
 \newtheorem{rem}[thm]{Remark}
 \newtheorem{note*}[thm]{Note}
  
 \numberwithin{equation}{section}

\providecommand{\Real}{\mathop{\rm Re}\nolimits}%
\providecommand{\clos}{\mathop{\rm clos}\nolimits}%
\providecommand{\id}{\mathop{\rm id}\nolimits}%

\usepackage{amssymb}

\begin{document}
% TOP MATTER

\title[Multivariable $\rho$-contractions]%
 {Multivariable $\rho$-contractions} % This is the full title of the paper

\author[Dmitry S. Kalyuzhny\u{\i}-Verbovetzki\u{\i}]{Dmitry S. Kalyuzhny\u{\i}-Verbovetzki\u{\i}}

\address{%
   Department of Mathematics\\
   The Weizmann Institute of Science\\
   Rehovot 76100\\
   Israel}

\email{dmitryk@wisdom.weizmann.ac.il}

\dedicatory{Dedicated to Israel Gohberg on his 75th birthday}

%Insert `2000 Mathematics Subject Classification' numbers here!


\begin{abstract}
We suggest a new version of the notion of $\rho$-dilation
($\rho>0$) of an $N$-tuple $\mathbf{A}=(A_1,\ldots,A_N)$ of
bounded linear operators on a common Hilbert space. We say that
$\mathbf{A}$ belongs to the class $C_{\rho,N}$ if $\mathbf{A}$
admits a $\rho$-dilation
$\widetilde{\mathbf{A}}=(\widetilde{A}_1,\ldots,\widetilde{A}_N)$
for which
$\zeta\widetilde{\mathbf{A}}:=\zeta_1\widetilde{A}_1+\cdots
+\zeta_N\widetilde{A}_N$ is a unitary operator for each
$\zeta:=(\zeta_1,\ldots,\zeta_N)$ in the unit torus
$\mathbb{T}^N$. For $N=1$ this class coincides with the class
$C_\rho$ of B. Sz.-Nagy and C. Foia\c{s}. We generalize the known
descriptions of $C_{\rho,1}=C_\rho$ to the case of $C_{\rho,N},\
N>1$, using so-called Agler kernels. Also, the notion of operator
radii $w_\rho, \rho>0$, is generalized to the case of $N$-tuples
of operators, and to the case of bounded (in a certain strong
sense) holomorphic operator-valued functions in the open unit
polydisk $\mathbb{D}^N$, with preservation of all the most
important their properties. Finally, we show that for each
$\rho>1$ and $N>1$ there exists an $\mathbf{A}=(A_1,\ldots,A_N)\in
C_{\rho,N}$ which is not simultaneously similar to any
$\mathbf{T}=(T_1,\ldots,T_N)\in C_{1,N}$, however if
$\mathbf{A}\in C_{\rho,N}$ admits a uniform unitary
$\rho$-dilation then $\mathbf{A}$ is simultaneously similar to
some $\mathbf{T}\in C_{1,N}$.
\end{abstract}

\maketitle

\section{Introduction}
\label{intro}

\noindent Linear pencils of operators $L_\mathbf{A}(z):=A_0+z_1A_1+\cdots +z_NA_N$
on a Hilbert space which take contractive (resp., unitary or
$J$-unitary for some signature operator $J=J^*=J^{-1}$) values for
all $z=(z_1,\ldots,z_N)$ in the \emph{unit torus}
$\mathbb{T}^N:=\{\zeta\in\mathbb{C}^N:|\zeta_k|=1,\
k=1,\ldots,N\}$ serve as one of possible generalizations of a
single contractive (resp., unitary, $J$-unitary) operator on a
Hilbert space. They appear in  constructions of Agler's unitary
colligation and corresponding conservative (unitary) scattering
$N$-dimensional discrete-time linear system of Roesser type
\cite{Ag,BT}, and also of Fornasini--Marchesini type \cite{BSV},
and dissipative (contractive), conservative (unitary) or
$J$-conservative ($J$-unitary) scattering $N$-dimensional linear
systems of one more form introduced in our paper \cite{K2} and
studied in \cite{K2,K3,K4,K5,K6,K7,BSV}. These constructions, in
particular, provide the transfer function realization formulae for
certain classes of holomorphic functions
\cite{Ag,BT,K2,K4,K7,BSV}, the solutions to the Nevanlinna--Pick
interpolation problem \cite{AgMc,BT}, the Toeplitz corona problem
\cite{AgMc,BT}, and the commutant lifting problem \cite{BLTT} in
several variables.

In \cite{K3} we developed the dilation theory for multidimensional
linear systems, and in particular gave a necessary and sufficient
condition for such a system to have a conservative dilation. As a
special case, this gave a criterion for the existence of a
unitary dilation of a contractive (on $\mathbb{T}^N$) linear
pencil of operators on a Hilbert space. Linear pencils of
operators satisfying this criterion inherit the most important
properties of single contraction operators on a Hilbert space
(note that, due to \cite{K1}, not all linear pencils which take
contractive operator values on $\mathbb{T}^N$ satisfy this
criterion).

The purpose of the present paper is to develop the theory of
$\rho$-contractions in several variables in the framework of
``linear pencils approach". We introduce the notion of
\emph{$\rho$-dilation} of an $N$-tuple
$\mathbf{A}=(A_1,\ldots,A_N)$ of bounded linear operators on a
common Hilbert space by means of a simultaneous $\rho$-dilation,
in the sense of B.~Sz.-Nagy and C.~Foia\c{s} \cite{SzNF1,SzNF}, of
the values of a homogeneous linear pencil of operators
$z\mathbf{A}:=\sum_{k=1}^Nz_kA_k$. The class $C_{\rho,N}$ consists
of those $N$-tuples of operators $\mathbf{A}=(A_1,\ldots,A_N)$
(\emph{$\rho$-contractions}) for which there exists a
$\rho$-dilation
$\widetilde{\mathbf{A}}=(\widetilde{A}_1,\ldots,\widetilde{A}_N)$
such that the operators
$\zeta\widetilde{\mathbf{A}}=\sum_{k=1}^N\zeta_k\widetilde{A}_k$
are unitary for all
$\zeta=(\zeta_1,\ldots,\zeta_N)\in\mathbb{T}^N$. On the one hand, this class
generalizes the class $C_{\rho,1}=C_\rho$ of
Sz.-Nagy and Foia\c{s} \cite{SzNF1,SzNF} consisting of operators
which admit a unitary $\rho$-dilation to the case $N>1$. On
the other hand, this class generalizes the class of $N$-tuples of operators $\mathbf{A}$
for which the associated linear pencil of operators $z\mathbf{A}$
admits a unitary dilation  in the sense of \cite{K3} (this
corresponds to $\rho=1$) to the case of $N$-tuples of operators
$\mathbf{A}$ which have a unitary $\rho$-dilation for $\rho\neq
1$.

The paper is organized as follows. Section~\ref{sec:prelim} gives
preliminaries on $\rho$-contracti\-ons for the case $N=1$. Namely,
we recall the relevant definitions, the known criteria for an
operator to be a $\rho$-contraction, i.e., to belong to the class
$C_{\rho}$ of Sz.-Nagy and Foia\c{s}, the notion of operator radii
$w_\rho$ and their properties, and the theorem on similarity of
$\rho$-contractions to contractions. In Section~\ref{sec:n-rho} we
give the definitions of a $\rho$-dilation of an $N$-tuple of
operators, and of the class $C_{\rho,N}$ of $\rho$-contractions for
the case $N>1$, and prove a theorem which generalizes the
criteria of $\rho$-contractiveness to this case, as well as to the
case $0<\rho\neq 1$. Some properties of classes $C_{\rho,N}$ are
discussed. Then it is shown that the notions of a $\rho$-contraction
and of the corresponding class $C_{\rho,N}$, as well as the 
theorem just mentioned, can be extended to holomorphic functions on the \emph{open unit polydisk}
$\mathbb{D}^N:=\{z\in\mathbb{C}^N:|z_k|<1,\ k=1,\ldots,N\}$ that are bounded in a certain
strong sense,
though the notion of unitary $\rho$-dilation is not relevant any
more in this case. In Section~\ref{sec:n-rad} we define operator
radii $w_{\rho,N}$ of $N$-tuples of operators, and
operator-function radii $w_{\rho,N}^{(\infty)}$ of bounded
holomorphic functions on $\mathbb{D}^N,\ \rho>0$. These radii
generalize $w_\rho$'s and inherit all the most important
properties of them. In Section~\ref{sec:similar} we prove that for
each $\rho>1$ and $N>1$ there exists an
$\mathbf{A}=(A_1,\ldots,A_N)\in C_{\rho,N}$ which is not
simultaneously similar to any $\mathbf{T}=(T_1,\ldots,T_N)\in
C_{1,N}$. Then we introduce the classes $C_{\rho,N}^u,\ \rho>0$,
of $N$-variable $\rho$-contractions $\mathbf{A}=(A_1,\ldots,A_N)$
which admit a uniform unitary $\rho$-dilation. We prove that if
$\mathbf{A}\in C_{\rho,N}^u$ for some $\rho>1$ then $\mathbf{A}$
is simultaneously similar to some $\mathbf{T}\in C_{1,N}^u$. Note,
that since the class $C_{\rho,N}^u$ (as well as $C_{\rho,N}$)
increases as a function of $\rho$, for any $\rho\leq 1$ an
$\mathbf{A}\in C_{\rho,N}^u$ (resp., $\mathbf{A}\in C_{\rho,N}$)
belongs to $C_{1,N}^u$ (resp., $C_{1,N}$) itself. We show the
relation of our results to ones of G.~Popescu \cite{Po2} where a
different notion of multivariable $\rho$-contractions has been
introduced, and the relevant theory has been developed. The
classes $C_{\rho,N}^u,\ \rho>0$, which appear in
Section~\ref{sec:similar} in connection with the similarity
problem discussed there, certainly deserve a further
investigation.

\section{Preliminaries}\label{sec:prelim}

Let $L(\mathcal{X,Y})$ denote the Banach space of bounded linear
operators mapping a  Hilbert space $\mathcal{X}$ into a Hilbert
space $\mathcal{Y}$, and $L(\mathcal{X}):=L(\mathcal{X,X})$. For
$\rho>0$, an operator $\widetilde{A}\in
L(\widetilde{\mathcal{X}})$ is said to be a \emph{$\rho$-dilation
of an operator} $A\in L(\mathcal{X})$ if
$\widetilde{\mathcal{X}}\supset\mathcal{X}$ and
\begin{equation}\label{eq:rhodil}
A^n=\rho P_\mathcal{X}\widetilde{A}^n|\mathcal{X},\quad
n\in\mathbb{N},
\end{equation}
 where $P_\mathcal{X}$ denotes the orthogonal projection onto the
subspace $\mathcal{X}$ in $\widetilde{\mathcal{X}}$. If, moreover,
$\widetilde{A}$ is a unitary operator then $\widetilde{A}$ is
called a \emph{unitary $\rho$-dilation of $A$}. In \cite{SzNF1}
(see also \cite{SzNF}) B.~Sz.-Nagy and C.~Foia\c{s} introduced the
classes $C_\rho,\ \rho>0$, consisting of operators which admit a
unitary $\rho$-dilation. Due to B.~Sz.-Nagy \cite{SzN}, the class
$C_1$ is precisely the class of all contractions, i.e., operators
$A$ such that $\| A\|\leq 1$. C.~A.~Berger \cite{Be} showed that
the class $C_2$ is precisely the class of all operators $A\in
L(\mathcal{X})$, for some Hilbert space $\mathcal{X}$, which have
the \emph{numerical radius}
\[ w(A)=\sup\{ |\langle Ax,x\rangle |:\ x\in\mathcal{X},\ \| x\|
=1\}\] equal to at most one. Thus, the classes $C_\rho,\ \rho>0$,
provide a framework for simultaneous investigation of these two
important classes of operators.

Recall that the \emph{Herglotz} (or \emph{Caratheodory})
\emph{class} $\mathcal{H(X)}$ (respectively, the \emph{Schur
class} $\mathcal{S(X)}$) consists of holomorphic functions $f$ on
the open unit disk $\mathbb{D}$ which take values in
$L(\mathcal{X})$ and satisfy $\Real f(z)=f(z)+f(z)^*\succeq 0$ in
the sense of positive semi-definiteness of an operator (resp., $\|
f(z)\|\leq 1)$ for all $z\in\mathbb{D}$. Let us recall some known
characterizations of the classes $C_\rho$.
\begin{thm}\label{thm:1-char}
Let $A\in L(\mathcal{X})$ and $\rho>0$. The following statements
are equivalent:
\begin{description}
   \item[(i)] $A\in C_\rho$;
    \item[(ii)] the function $k_\rho^A(z,w):=\rho
    I_\mathcal{X}-(\rho-1)\left((zA+(wA)^*\right)+(\rho-2)(wA)^*zA$ satisfies
$k_\rho^A(z,z)\succeq 0$ for all $z\in\clos(\mathbb{D})$;
    \item[(iii)] the function
    $\psi_\rho^A(z):=(1-\frac{2}{\rho})I_\mathcal{X}+\frac{2}{\rho}(I_\mathcal{X}-zA)^{-1}$
    belongs to $\mathcal{H(X)}$;
    \item[(iv)] the function
    $\varphi_\rho^A(z):=zA\left((\varrho-1)zA-\rho I_\mathcal{X}\right)^{-1}$
    belongs to $\mathcal{S(X)}$.
\end{description}
\end{thm}
Conditions (ii) and (iii) of Theorem~\ref{thm:1-char} each
characterizing the class $C_\rho$ appear in \cite{SzNF1}, while
condition (iv) is due to C.~Davis \cite{Da}.
\begin{cor}\label{cor:1-char'}
 Condition $(${\rm ii}$)$ in Theorem~\ref{thm:1-char} can be replaced
by \vspace{2mm}
\begin{description}
    \item[(ii')] $k_\rho^A(C,C):=\rho I_\mathcal{X}\otimes
    I_{\mathcal{H}_C}-(\rho-1)(A\otimes C+
    (A\otimes
    C)^*)\\+(\rho-2)(A\otimes C)^*(A\otimes C)\succeq 0$\\ for any
    contraction $C$ on a Hilbert space $\mathcal{H}_C$.
\end{description}
\end{cor}
\begin{proof} Indeed, (ii')$\Rightarrow$(ii), hence
(ii')$\Rightarrow$(i). Conversely, if $A\in C_\rho\cap
L(\mathcal{X})$ then for any contraction $C$ on $\mathcal{H}_C$
one has $A\otimes C\in C_\rho$ because, by \cite{SzN}, $C$ admits
a unitary dilation $\widetilde{C}$, and $A$ admits a unitary
$\rho$-dilation $\widetilde{A}$, thus
$\widetilde{A}\otimes\widetilde{C}$ is a unitary $\rho$-dilation
of $A\otimes C$:
\begin{eqnarray*}
(A\otimes C)^n &=& A^n\otimes C^n= (\rho
P_\mathcal{X}\widetilde{A}^n|\mathcal{X})\otimes
(P_{\mathcal{H}_C}\widetilde{C}^n|\mathcal{H}_C)\\
&=& \rho P_{\mathcal{X\otimes
H}_C}(\widetilde{A}^n\otimes\widetilde{C}^n)|\mathcal{X\otimes
H}_C\\
 &=& \rho P_{\mathcal{X\otimes
H}_C}(\widetilde{A}\otimes\widetilde{C})^n|\mathcal{X\otimes
H}_C,\quad n\in\mathbb{N}.
\end{eqnarray*}
Therefore, $k_\rho^A(C,C)=k_\rho^{A\otimes C}(1,1)\succeq 0$,
i.e., (ii') is valid.
\end{proof}
\begin{cor}\label{cor:ac}
 Condition \vspace{2mm}
\begin{description}
    \item[(v)] $A\otimes C\in C_\rho$ for any
    contraction $C$ on a Hilbert space,
\end{description}\vspace{2mm}
 is equivalent to each of conditions $(${\rm i}$)$--$(${\rm iv}$)$ of
Theorem~\ref{thm:1-char}.
\end{cor}
\begin{proof}
See the proof of Corollary~\ref{cor:1-char'}.
\end{proof}
 Any operator $A\in C_\rho$ is
\emph{power-bounded}:
\begin{equation}\label{eq:1-pb}
 \| A^n\|\leq \rho,\quad n\in\mathbb{N},
\end{equation}
moreover, its \emph{spectral radius}
\begin{equation}\label{eq:1-sp-rad}
\nu(A)=\lim_{n\rightarrow+\infty}\| A^n\|^{\frac{1}{n}}
\end{equation}
 is at most
one. In \cite{SzNF1} an example of a power-bounded operator which
is not contained in any of the classes $C_\rho,\ \rho>0$, is
given. However, J.~A.~R.~Holbrook \cite{H1} showed that any
bounded linear operator $A$ with $\nu(A)\leq 1$ can be
approximated in the operator norm topology by elements of the
classes $C_\rho$. More precisely, if $C_\infty$ denotes the class
of bounded linear operators with spectral radius at most one, and
$\mathcal{X}$ is a Hilbert space, then
\begin{equation}\label{eq:1-clos}
 C_\infty\cap
L(\mathcal{X})=\clos\left\{\bigcup_{0<\rho<\infty}C_\rho\cap
L(\mathcal{X})\right\}.
\end{equation}

For a fixed Hilbert space $\mathcal{X}$, the class $C_\rho$ as a
function of $\rho$ increases \cite{SzNF1}:
\begin{equation}\label{eq:incr}
C_\rho\subset C_{\rho'}\ \mbox{for}\ \rho<\rho'.
\end{equation}
Moreover, it was shown by E.~Durszt \cite{Du} that $C_\rho$
increases strictly for $\dim\mathcal{X}\geq 2$:
\[ C_\rho\neq C_{\rho'}\ \mbox{for}\ \rho\neq\rho'.\]
\begin{prop}\label{prop:1c}
For $\mathcal{X}=\mathbb{C}$, the classes $C_\rho$ coincide for
all $\rho\geq 1$, and strictly increase for $0<\rho<1$:
\[ C_\rho\varsubsetneq C_{\rho'}\ \mbox{for}\ 0<\rho<\rho'\leq
1.\]
\end{prop}
\begin{proof}
If $a\in\mathbb{C}\cong L(\mathbb{C})$ belongs to $C_\rho$ then
$\|a\|=|a| =\nu(a)\leq 1$. Hence $C_\rho\subset C_1$ for any
$\rho>0$. Since (\ref{eq:incr}) implies $C_\rho\supset C_1$ for
$\rho\geq 1$, we get $C_\rho=C_1$ for this case, that proves the
first part of this proposition.

For the proof of the second part, we will show that for any
$\varepsilon,\rho:0<\varepsilon<\rho<1$, one has
\begin{equation}\label{eq:a}
 a:=\frac{\rho}{2-\rho}\in C_\rho\backslash
C_{\rho-\varepsilon}.
\end{equation}
 If $0\leq\varepsilon<\rho$ then, by condition (ii) in
Theorem~\ref{thm:1-char}, the inclusion $a\in
C_{\rho-\varepsilon}$ is equivalent to
\[
\rho-\varepsilon-(\rho-\varepsilon-1)(az+\bar{a}\bar{z})+(\rho-\varepsilon-2)|az|^2\geq
0,\quad z\in\clos(\mathbb{D}),\] which for $a=\frac{\rho}{2-\rho}$
turns into
\[\rho-\varepsilon-2(\rho-\varepsilon-1)\frac{\rho}{2-\rho}r\cos\theta+(\rho-\varepsilon-2)\left(\frac{\rho r}{2-\rho}\right)^2\geq
0,\quad r\in [0,1],\theta\in [0,2\pi).\] Since
$\rho-\varepsilon-1<0$, the left-hand side of this inequality, as
a function of $\theta$ for a fixed $r$, has a minimum at
$\theta=\pi$, so the latter condition turns into
\[\rho-\varepsilon+2(\rho-\varepsilon-1)\frac{\rho r}{2-\rho}+(\rho-\varepsilon-2)\left(\frac{\rho r}{2-\rho}\right)^2\geq
0,\quad r\in [0,1].\] The left-hand side attains its minimum at
$r=1$, thus the latter inequality turns into
\[\rho-\varepsilon+2(\rho-\varepsilon-1)\frac{\rho}{2-\rho}+(\rho-\varepsilon-2)\left(\frac{\rho}{2-\rho}\right)^2=-\frac{4\varepsilon}{(2-\rho)^2}\geq
0,\] which is possible if and only if $\varepsilon=0$. Thus,
(\ref{eq:a}) is true.
\end{proof}
The properties of the classes $C_\rho$ become more clear due to
the following numerical characteristics of operators.
J.~A.~R.~Holbrook \cite{H1} and J.~P.~Williams \cite{W},
independently, introduced for any $A\in L(\mathcal{X})$ the
\emph{operator radii}
\begin{equation}\label{eq:1-rad}
 w_\rho(A):=\inf\{ u>0:\frac{1}{u}A\in C_\rho\}.
\end{equation}
\begin{thm}\label{thm:1-rad}
$w_\rho(\cdot)$ has the following properties:
\begin{description}
    \item[(i)] $w_\rho(A)<\infty$;
    \item[(ii)] $w_\rho(A)>0$ unless $A=0$, moreover,
    $w_\rho(A)\geq\frac{1}{\rho}\| A\|$;
    \item[(iii)] $\forall\mu\in\mathbb{C},\quad w_\rho(\mu
    A)=|\mu|w_\rho(A)$;
    \item[(iv)] $w_\rho(A)\leq 1$ if and only if $A\in C_\rho$;
    \item[(v)] $w_\rho(\cdot)$ is a norm on $L(\mathcal{X})$ for any $\rho:\ 0<\rho\leq 2$,
    and not a norm on $L(\mathcal{X}),\ \dim\mathcal{X}\geq 2$, for any
    $\rho>2$;
    \item[(vi)] $w_1(A)=\| A\|$ (of course, here $\|\cdot\|$ is the
    operator norm on $L(\mathcal{X})$ with respect to the
    Hilbert-space metric on $\mathcal{X}$);
    \item[(vii)] $w_2(A)=w(A)$;
    \item[(viii)] $w_\infty(A):=\lim\limits_{\rho\rightarrow
    +\infty}w_\rho(A)=\nu(A)$;
    \item[(ix)] $w_\rho(I_\mathcal{X})=\left\{\begin{array}{ll}
      1 & \mbox{for}\ \rho\geq 1, \\
      \frac{2}{\rho}-1 & \mbox{for}\ 0<\rho<1;
    \end{array}\right.$
    \item[(x)] if $0<\rho<\rho'$ then $w_{\rho'}(A)\leq
    w_\rho(A)\leq\left(\frac{2\rho'}{\rho}-1\right)w_{\rho'}(A)$,
    thus $w_\rho(A)$ is continuous in $\rho$ and non-increasing as
    $\rho$ increases;
    \item[(xi)] if $\| A\|=1$ and $A^2=0$ then, for any $\rho>0$,
    $w_\rho(A)=\frac{1}{\rho}$;
    \item[(xii)] if for some $\rho_0$ one has $w_{\rho_0}(A)>w_\infty(A)\
    (=\nu(A))$ then for any $\rho>\rho_0$ one has
    $w_{\rho_0}(A)>w_\rho(A)$;
    \item[(xiii)] $\lg w_\rho(A)$ is a convex function in $\rho,\
    0<\rho<+\infty$;
    \item[(xiv)] $w_\rho(A)$ is a convex function in $\rho,\
    0<\rho<+\infty$;
    \item[(xv)] the function $h_A(\rho):=\rho w_\rho(A)$ is
    non-decreasing on $[1,+\infty)$, and non-increa\-sing on
    $(0,1)$;
    \item[(xvi)] for any $\rho$ such that $0<\rho<2$ one has $\rho
    w_\rho(A)=(2-\rho)w_{2-\rho}(A)$, and $\lim\limits_{\rho\downarrow
    0}\frac{\rho}{2}w_\rho(A)=w_2(A)\ (=w(A))$;
    \item[(xvii)] $\forall\rho:\ 0<\rho\leq 1,\quad
    w_\rho(A)\geq\left(\frac{2}{\rho}-1\right)w_2(A)$;
    \item[(xviii)] $\forall A,B\in L(\mathcal{X}),\ \forall\rho\geq 1,\quad
    w_\rho(AB)\leq\rho^2w_\rho(A)w_\rho(B)$, moreover, $\rho^2$ is
    the best constant in this inequality for the case $\dim\mathcal{X}\geq
    2$;
    \item[(xix)] $\forall A,B\in L(\mathcal{X}),\ \forall\rho:\ 0<\rho<1,\quad
    w_\rho(AB)\leq(2-\rho)\rho w_\rho(A)w_\rho(B)$, moreover, $(2-\rho)\rho$ is
    the best constant in this inequality for the case $\dim\mathcal{X}\geq
    2$;
    \item[(xx)] $\forall\rho>0,\ \forall n\in\mathbb{N},\quad
    w_\rho(A^n)\leq w_\rho(A)^n$.
\end{description}
\end{thm}
Properties (i)--(xii), (xviii), and (xx) were proved by
J.~A.~R.~Holbrook \cite{H1}, properties (xiii)--(xvi) were
discovered by T.~Ando and K.~Nishio \cite{AnN}. Property (xix) was
shown by K.~Okubo and T.~Ando \cite{OAn}, and follows also from
(xvi) and (xviii). Finally, property (xvii) easily follows from
(x) and (xvi). Indeed, for $0<\rho\leq 1$ one has
$w_{2-\rho}(A)\geq w_2(A)$, hence $\rho
w_\rho(A)=(2-\rho)w_{2-\rho}(A)\geq(2-\rho)w_2(A)$, which implies
(xvii).

We have listed in Theorem~\ref{thm:1-rad} only the most important,
as it seems to us, properties of operator radii $w_\rho(\cdot)$.
Other properties of $w_\rho(\cdot)$ can be found in
\cite{H1,H2,FH,AnN,OAn,BaCa} and elsewhere.

Let us note that properties of the classes $C_\rho$ discussed
before Theorem~\ref{thm:1-rad}, including
Proposition~\ref{prop:1c}, can be deduced from properties (iv),
(vi)--(x) in Theorem~\ref{thm:1-rad}. Due to property (iv) in
Theorem~\ref{thm:1-rad}, operators from the classes $C_\rho$ are
called \emph{$\rho$-contractions}.

Any $A\in C_\rho$ satisfies the following \emph{generalized von
Neumann inequality} \cite{SzNF1}: for any polynomial $p$ of one
variable
\begin{equation}\label{eq:vn}
\| p(A)\|\leq\max_{|z|\leq 1}|\rho p(z)+(1-\rho)p(0)|.
\end{equation}
Let $A\in L(\mathcal{X}), B\in L(\mathcal{Y})$. Then $A$ is said
to be \emph{similar} to $B$ if there exists a bounded invertible
operator $S\in L(\mathcal{X,Y})$ such that
\begin{equation}\label{eq:1-similar}
A=S^{-1}BS.
\end{equation}
 B.~Sz.-Nagy and C.~Foia\c{s} proved in \cite{SzNF2}
(see also \cite{SzNF}) that any $A\in C_\rho$ is similar to some
$T\in C_1$, i.e., any $\rho$-contraction is similar to a
contraction.

To conclude this section, let us remark that the classes
$C_\rho$ are of continuous interest, e.g., see recent works
\cite{DMcCW,CaF,BaCa,NO,OS}. In \cite{Po2} the classes $C_\rho$
were extended to a multivariable setting; we shall discuss this generalization in
Section~\ref{sec:similar}.

\section{The classes $C_{\rho,N}$}\label{sec:n-rho}
Let $\rho>0$. We will say that an $N$-tuple of operators
$\widetilde{\mathbf{A}}=(\widetilde{A}_1,\ldots,\widetilde{A}_N)\in
L(\widetilde{\mathcal{X}})^N$ is a \emph{$\rho$-dilation of an
$N$-tuple of operators} $\mathbf{A}=(A_1,\ldots,A_N)\in
L(\mathcal{X})^N$ if $\widetilde{\mathcal{X}}\supset\mathcal{X}$,
and for any $z=(z_1,\ldots,z_N)\in\mathbb{C}^N$ the operator
$z\widetilde{\mathbf{A}}=\sum_{k=1}^Nz_k\widetilde{A}_k$ is a
$\rho$-dilation, in the sense of \cite{SzNF1}, of the operator
$z\mathbf{A}=\sum_{k=1}^Nz_kA_k$, i.e.,
\begin{equation}\label{eq:n-rhodil}
 (z\mathbf{A})^n=\rho
P_\mathcal{X}(z\widetilde{\mathbf{A}})^n|\mathcal{X},\quad
z\in\mathbb{C}^N,\ n\in\mathbb{N}.
\end{equation}
 These relations are equivalent to
\begin{equation}\label{eq:n-rhodil'}
 \mathbf{A}^t=\rho
P_\mathcal{X}\widetilde{\mathbf{A}}^t|\mathcal{X},\quad
t\in\mathbb{Z}^N_+:=\{\tau\in\mathbb{Z}^N:\tau_k\geq 0,\
k=1,\ldots,N\},
\end{equation}
 where $\mathbf{A}^t,\ t\in\mathbb{Z}^N_+$, are \emph{symmetrized
multi-powers of} $\mathbf{A}$:
\[\mathbf{A}^t:=\frac{t!}{|t|!}\sum_\sigma A_{[\sigma(1)]}\cdots
A_{[\sigma(|t|)]},\] and analogously for $\widetilde{\mathbf{A}}$.
Here for a multi-index $t=(t_1,\ldots,t_N),\
t!:=t_1!\cdots t_N!$ and $|t|:=t_1+\cdots+t_N$;
$\sigma$ runs over the set of all permutations with repetitions in a string of $|t|$ numbers from the set $\{ 1,\ldots,N\}$ such that the
$\kappa$-th number $[\kappa]\in\{1,\ldots,N\}$ appears in this string
$t_{[\kappa]}$ times. Say, if $t=(1,2,0,\ldots,0)$ then
\[ \mathbf{A}^t=\frac{A_1A_2^2+A_2A_1A_2+A_2^2A_1}{3}.\]
In the case of a commutative $N$-tuple $\mathbf{A}$ one has
$\mathbf{A}^t=A_1^{t_1}\cdots A_N^{t_N}$, i.e., a
usual multi-power. \begin{note*} Compare (\ref{eq:n-rhodil}) and
(\ref{eq:n-rhodil'}) with (\ref{eq:rhodil}).
\end{note*}
In the case $\rho=1$ the notion of $\rho$-dilation of an $N$-tuple
of operators $\mathbf{A}=(A_1,\ldots,A_N)$ coincides with the
notion of dilation of $\mathbf{A}$ (or corresponding linear pencil
$z\mathbf{A}$) as defined in \cite{K3}.

We will call $\widetilde{\mathbf{A}}\in
L(\widetilde{\mathcal{X}})^N$ a \emph{unitary $\rho$-dilation of
$\mathbf{A}\in L(\mathcal{X})^N$} if $\widetilde{\mathbf{A}}$ is a
$\rho$-dilation of $\mathbf{A}$ and for any $\zeta\in\mathbb{T}^N$
the operator
$\zeta\widetilde{\mathbf{A}}=\sum_{k=1}^N\zeta_k\widetilde{A}_k$
is unitary. The class of operator $N$-tuples which admit a unitary
$\rho$-dilation will be denoted by $C_{\rho,N}$.

Let $\mathcal{C}^N$ denote the family of all $N$-tuples
$\mathbf{C}=(C_1,\ldots,C_N)$ of \emph{commuting strict
contractions} on a common Hilbert space $\mathcal{H}_\mathbf{C}$,
i.e., $C_kC_j=C_jC_k$ and $\| C_k\|<1$ for all $k,j\in\{
1,\ldots,N\}$. An $L(\mathcal{X})$-valued function
\[
k(z,w)=\sum_{(t,s)\in\mathbb{Z}^N_+\times\mathbb{Z}^N_+}\hat{k}(t,s)\bar{w}^sz^t,\quad
(z,w)\in\mathbb{D}^N\times\mathbb{D}^N,\] which is holomorphic in
$z\in\mathbb{D}^N$ and anti-holomorphic in $w\in\mathbb{D}^N$,
will be called an \emph{Agler kernel} if
\begin{equation}\label{eq:n-calc}
k(\mathbf{C},\mathbf{C}):=\sum_{(t,s)\in\mathbb{Z}^N_+\times\mathbb{Z}^N_+}
\hat{k}(t,s)\otimes\mathbf{C}^{*s}\mathbf{C}^t\succeq 0,\quad
\mathbf{C}\in\mathcal{C}^N,
\end{equation}
  where the series converges in the
operator norm topology on $L(\mathcal{X\otimes H}_\mathbf{C})$.
The \emph{Agler--Herglotz class $\mathcal{AH}_N(\mathcal{X})$}
(resp., the \emph{Agler--Schur class
$\mathcal{AS}_N(\mathcal{X})$}) is the class of all
$L(\mathcal{X})$-valued functions $f$ holomorphic on
$\mathbb{D}^N$ for which $k(z,w)=f(z)+f(w)^*$ (resp.,
$k(z,w)=I_\mathcal{X}-f(w)^*f(z)$) is an Agler kernel. Agler
kernels, as well as the classes $\mathcal{AH}_N(\mathcal{X})$ and
$\mathcal{AS}_N(\mathcal{X})$, were defined and studied by
J.~Agler in \cite{Ag}. The von Neumann inequality \cite{vN}
implies that $\mathcal{AS}_1(\mathcal{X})=\mathcal{S(X)}$ and
$\mathcal{AH}_1(\mathcal{X})=\mathcal{H(X)}$.
\begin{rem}\label{rem:a-ker}
 The function $k_\rho^A(z,w)$ from condition (ii) in
Theorem~\ref{thm:1-char}, due to Corollary~\ref{cor:1-char'}, is
an Agler kernel ($N=1$).
\end{rem}
\begin{thm}\label{thm:n-char}
Let $\mathbf{A}\in L(\mathcal{X})^N,\ \rho>0$. The following
conditions are equivalent:
\begin{description}
    \item[(i)] $\mathbf{A}\in C_{\rho,N}$;
    \item[(ii)] the function $k_{\rho,N}^\mathbf{A}(z,w):=\rho
I_\mathcal{X}-(\rho-1)\left((z\mathbf{A}+(w\mathbf{A})^*\right)+(\rho-2)(w\mathbf{A})^*z\mathbf{A}$
is an Agler kernel on $\mathbb{D}^N\times\mathbb{D}^N$;
    \item[(iii)] the function
$\psi_{\rho,N}^\mathbf{A}(z):=(1-\frac{2}{\rho})I_\mathcal{X}+\frac{2}{\rho}(I_\mathcal{X}-z\mathbf{A})^{-1}$
    belongs to $\mathcal{AH}_N(\mathcal{X})$;
    \item[(iv)] the function
    $\varphi_{\rho,N}^\mathbf{A}(z):=z\mathbf{A}\left((\varrho-1)z\mathbf{A}-\rho I_\mathcal{X}\right)^{-1}$
    belongs to $\mathcal{AS}_N(\mathcal{X})$;
    \item[(v)] $\mathbf{A\otimes C}:=\sum_{k=1}^NA_k\otimes C_k\in
    C_\rho=C_{\rho,1}$ for all $\mathbf{C}\in\mathcal{C}^N.$
\end{description}
\end{thm}
\begin{rem}\label{rem:general}
This theorem generalizes Theorem~\ref{thm:1-char} with condition
(ii) replaced by condition (ii') from Corollary~\ref{cor:1-char'},
and added condition (v) from Corollary~\ref{cor:ac}.
\end{rem}
\begin{proof}[Proof of Theorem~\ref{thm:n-char}]
(i)$\Leftrightarrow$(iii). The proof of this part combines the
idea of B.~Sz.-Nagy and C.~Foia\c{s} \cite{SzNF1} for the proof of
the equivalence (i)$\Leftrightarrow$(iii) in
Theorem~\ref{thm:1-char} (see Remark~\ref{rem:general}) with
Agler's representation of functions from
$\mathcal{AH}_N(\mathcal{X})$ \cite{Ag}. Let
$\mathbf{A}=(A_1,\ldots,A_N)\in C_{\rho,N}\cap L(\mathcal{X})^N$,
and
$\widetilde{\mathbf{A}}=(\widetilde{A}_1,\ldots,\widetilde{A}_N)\in
L(\widetilde{\mathcal{X}})^N$ be a unitary $\rho$-dilation of
$\mathbf{A}$. By Corollary~4.3 in \cite{K3}, the linear function
$L_{\widetilde{\mathbf{A}}}(z)=z\widetilde{\mathbf{A}}$ belongs to
the class $\mathcal{AS}_N(\widetilde{\mathcal{X}})$. Since for any
$\mathbf{C}\in\mathcal{C}^N$ one has
$(1+\varepsilon)\mathbf{C}\in\mathcal{C}^N$ for a sufficiently
small $\varepsilon>0$, the operator
$\widetilde{\mathbf{A}}\otimes\mathbf{C}$, as well as
$\widetilde{\mathbf{A}}\otimes(1+\varepsilon)\mathbf{C}$, is
contractive. Thus, $\widetilde{\mathbf{A}}\otimes\mathbf{C}$ is a
strict contraction, and the series
\[
I_{\widetilde{\mathcal{X}}\otimes\mathcal{H}_\mathbf{C}}+2\sum_{n=1}^\infty(\widetilde{\mathbf{A}}\otimes\mathbf{C})^n\]
converges in the
$L(\widetilde{\mathcal{X}}\otimes\mathcal{H}_\mathbf{C})$-norm to
\[
(I_{\widetilde{\mathcal{X}}\otimes\mathcal{H}_\mathbf{C}}+\widetilde{\mathbf{A}}\otimes\mathbf{C})
(I_{\widetilde{\mathcal{X}}\otimes\mathcal{H}_\mathbf{C}}-\widetilde{\mathbf{A}}\otimes\mathbf{C})^{-1}.\]
Moreover,
\begin{equation}\label{eq:pos-real}
    \Real[(I_{\widetilde{\mathcal{X}}\otimes\mathcal{H}_\mathbf{C}}+\widetilde{\mathbf{A}}\otimes\mathbf{C})
(I_{\widetilde{\mathcal{X}}\otimes\mathcal{H}_\mathbf{C}}-\widetilde{\mathbf{A}}\otimes\mathbf{C})^{-1}]\succeq
0.
\end{equation}
Therefore,
\begin{eqnarray*}
\lefteqn{P_{\mathcal{X\otimes
H}_\mathbf{C}}(I_{\widetilde{\mathcal{X}}\otimes\mathcal{H}_\mathbf{C}}+\widetilde{\mathbf{A}}\otimes\mathbf{C})
(I_{\widetilde{\mathcal{X}}\otimes\mathcal{H}_\mathbf{C}}-\widetilde{\mathbf{A}}\otimes\mathbf{C})^{-1}|\mathcal{X\otimes
H}_\mathbf{C} }\\
&=&\left. P_{\mathcal{X\otimes
H}_\mathbf{C}}\left(I_{\widetilde{\mathcal{X}}\otimes\mathcal{H}_\mathbf{C}}+2\sum_{n=1}^\infty(\widetilde{\mathbf{A}}\otimes\mathbf{C})^n\right)\right|
\mathcal{X\otimes
H}_\mathbf{C}\\
&=&
I_{\mathcal{X}\otimes\mathcal{H}_\mathbf{C}}+2\sum_{n=1}^\infty\sum_{|t|=n}\frac{n!}{t!}(P_\mathcal{X}\otimes
I_{\mathcal{H}_\mathbf{C}})(\widetilde{\mathbf{A}}^t\otimes\mathbf{C}^t)|\mathcal{X\otimes
H}_\mathbf{C}\\
 &=&
I_{\mathcal{X}\otimes\mathcal{H}_\mathbf{C}}+\frac{2}{\rho}\sum_{n=1}^\infty\sum_{|t|=n}\frac{n!}{t!}\mathbf{A}^t\otimes\mathbf{C}^t
=
I_{\mathcal{X}\otimes\mathcal{H}_\mathbf{C}}+\frac{2}{\rho}\sum_{n=1}^\infty(\mathbf{A\otimes
C})^n \\
&=&
(1-\frac{2}{\rho})I_{\mathcal{X}\otimes\mathcal{H}}+\frac{2}{\rho}(I_{\mathcal{X}\otimes\mathcal{H}}-\mathbf{A\otimes
C})^{-1}=\psi_{\rho,N}^\mathbf{A}(\mathbf{C}),
\end{eqnarray*}
and (\ref{eq:pos-real}) implies
$\Real\psi_{\rho,N}^\mathbf{A}(\mathbf{C})\succeq 0$. Since the
function
$(I_{\widetilde{\mathcal{X}}}+z\widetilde{\mathbf{A}})(I_{\widetilde{\mathcal{X}}}-z\widetilde{\mathbf{A}})^{-1}$
is well-defined and holomorphic on $\mathbb{\mathbb{D}}^N$, so is
\begin{equation}\label{eq:psi}
\psi_{\rho,N}^\mathbf{A}(z)=P_\mathcal{X}(I_{\widetilde{\mathcal{X}}}+
z\widetilde{\mathbf{A}})(I_{\widetilde{\mathcal{X}}}-z\widetilde{\mathbf{A}})^{-1}|\mathcal{X},\quad
z\in\mathbb{D}^N,
\end{equation}
and we obtain
$\psi_{\rho,N}^\mathbf{A}\in\mathcal{AH}_N(\mathcal{X})$.

Conversely, let
$\psi_{\rho,N}^\mathbf{A}\in\mathcal{AH}_N(\mathcal{X})$. Since
$\psi_{\rho,N}^\mathbf{A}(0)=I_\mathcal{X}$, according to
\cite{Ag}, there exist a Hilbert space
$\widetilde{\mathcal{X}}\supset\mathcal{X}$, its subspaces
$\widetilde{\mathcal{X}}_1,\ldots,\widetilde{\mathcal{X}}_N$
satisfying
$\widetilde{\mathcal{X}}=\bigoplus_{k=1}^N\widetilde{\mathcal{X}}_k$,
and a unitary operator $U\in L(\widetilde{\mathcal{X}})$ such that
\begin{equation}\label{eq:ah}
\psi_{\rho,N}^\mathbf{A}(z)=P_\mathcal{X}(I_{\widetilde{\mathcal{X}}}+U(z\mathbf{P}))(I_{\widetilde{\mathcal{X}}}-U(z\mathbf{P}))^{-1}|\mathcal{X},\quad
z\in \mathbb{D}^N,
\end{equation}
 where
$z\mathbf{P}:=\sum_{k=1}^Nz_kP_{\widetilde{\mathcal{X}}_k}$, i.e.,
we get (\ref{eq:psi}) with
$\widetilde{A}_k=UP_{\widetilde{\mathcal{X}}_k},\ k=1,\ldots,N$.
Note that for each $\zeta\in\mathbb{T}^N$ the operator
$\zeta\widetilde{\mathbf{A}}$ is unitary. Developing both parts of
(\ref{eq:ah}) into the series in homogeneous polynomials
convergent in the operator norm, we get
\[
I_\mathcal{X}+\frac{2}{\rho}\sum_{n=1}^\infty(z\mathbf{A})^n=I_\mathcal{X}+2\sum_{n=1}^\infty
P_\mathcal{X}(z\widetilde{\mathbf{A}})^n|\mathcal{X},\quad
z\in\mathbb{D}^N,\] that implies the relations
\[ (z\mathbf{A})^n=\rho
P_\mathcal{X}(z\widetilde{\mathbf{A}})^n|\mathcal{X},\quad
n\in\mathbb{N},\] for all $z\in\mathbb{D}^N$, and hence for all
$z\in\mathbb{C}^N$. Thus, $\widetilde{\mathbf{A}}$ is a unitary
$\rho$-dilation of $\mathbf{A}$, and $\mathbf{A}\in C_{\rho,N}$.
The equivalence (i)$\Leftrightarrow$(iii) is proved.

 Note that in
this proof we have established that each Agler representation
(\ref{eq:ah}) of $\psi_{\rho,N}^\mathbf{A}$ gives rise to a
unitary $\rho$-dilation $\widetilde{\mathbf{A}}$ of $\mathbf{A}$,
and vice versa. Indeed, we already showed that (\ref{eq:ah})
determines $\widetilde{\mathbf{A}}$. Conversely, if
$\widetilde{\mathbf{A}}\in L(\widetilde{\mathcal{X}})^N$ is a
unitary $\rho$-dilation of $\mathbf{A}$, then (\ref{eq:psi})
holds. Set $U:=\sum_{k=1}^N\widetilde{A}_k\in
L(\widetilde{\mathcal{X}})$ and
$\widetilde{\mathcal{X}}_k:=\widetilde{A}_k^*\widetilde{\mathcal{X}},\
k=1,\ldots,N$. Then $U$ is unitary, $\widetilde{\mathcal{X}}_k$ is
a closed subspace in $\widetilde{\mathcal{X}}$ for each
$k=1,\ldots,N$, the subspaces $\widetilde{\mathcal{X}}_k$ are
pairwise orthogonal, and
$\widetilde{\mathcal{X}}=\bigoplus_{k=1}^N\widetilde{\mathcal{X}}_k$
(see Proposition~2.4 in \cite{K2}). Thus, (\ref{eq:psi}) turns
into (\ref{eq:ah}).

(v)$\Leftrightarrow$(iv). Let (v) be true. By
Theorem~\ref{thm:1-char} applied for $\mathbf{A\otimes C}$ with a
$\mathbf{C}\in\mathcal{C}^N$, one has
$\varphi_{\rho}^{\mathbf{A\otimes
C}}\in\mathcal{S}(\mathcal{X\otimes H}_\mathbf{C})$. For
$\varepsilon>0$ small enough,
$(1+\varepsilon)\mathbf{C}\in\mathcal{C}^N$, hence
$\mathbf{A}\otimes(1+\varepsilon)\mathbf{C}\in C_\rho$, and
$\varphi_\rho^{\mathbf{A}\otimes(1+\varepsilon)\mathbf{C}}\in\mathcal{S}(\mathcal{X\otimes
H}_\mathbf{C})$. Thus,
\[ \varphi_{\rho,N}^\mathbf{A}(\mathbf{C})=\mathbf{A\otimes
C}((\rho-1)\mathbf{A\otimes C}-\rho
I_{\mathcal{X}\otimes\mathcal{H}_\mathbf{C}})^{-1}=\varphi_\rho^{\mathbf{A}\otimes(1+\varepsilon)\mathbf{C}}\left(\frac{1}{1+\varepsilon}\right)\]
is a contraction on $\mathcal{X}\otimes\mathcal{H}_\mathbf{C}$. In
particular, $\varphi_{\rho,N}^\mathbf{A}(z)$ is well-defined,
holomorphic and contractive on $\mathbb{D}^N$. Finally,
$\varphi_{\rho,N}^\mathbf{A}\in\mathcal{AS}_N(\mathcal{X})$.

Conversely, if (iv) is true then for any
$\mathbf{C}\in\mathcal{C}^N$:
\[\varphi_\rho^{\mathbf{A\otimes C}}(\lambda)=\lambda\mathbf{A\otimes
C}((\rho-1)\lambda\mathbf{A\otimes C}-\rho
I_{\mathcal{X}\otimes\mathcal{H}_\mathbf{C}})^{-1}=\varphi_{\rho,N}^\mathbf{A}(\lambda\mathbf{C})\]
is well-defined, holomorphic and contractive for
$\lambda\in\mathbb{D}$. Thus, $\varphi_\rho^{\mathbf{A\otimes
C}}\in\mathcal{S}(\mathcal{X\otimes H}_\mathbf{C})$, and by
Theorem~\ref{thm:1-char}, $\mathbf{A\otimes C}\in C_\rho$.

(v)$\Leftrightarrow$(iii) and (v)$\Leftrightarrow$(ii) are proved
analogously, using the following relations for
$\mathbf{C}\in\mathcal{C}^N,\ \lambda\in\mathbb{D}$:
\begin{eqnarray*}
\psi_{\rho,N}^\mathbf{A}(\mathbf{C})=\psi_\rho^{\mathbf{A}\otimes(1+\varepsilon)\mathbf{C}}\left(\frac{1}{1+\varepsilon}\right),
\quad\psi_\rho^{\mathbf{A\otimes
C}}(\lambda)=\psi_{\rho,N}^\mathbf{A}(\lambda\mathbf{C}),\\
k_{\rho,N}^\mathbf{A}(\mathbf{C,C})=k_\rho^{\mathbf{A}\otimes\mathbf{C}}(1,1),
\quad k_\rho^{\mathbf{A\otimes
C}}(\lambda,\lambda)=k_{\rho,N}^\mathbf{A}(\lambda\mathbf{C},\lambda\mathbf{C}).
\end{eqnarray*}
The proof is complete.
\end{proof}
\begin{rem}\label{rem:1-contr}
For the case $\rho=1$ each of conditions (ii)--(v) in
Theorem~\ref{thm:n-char} means that for any
$\mathbf{C}\in\mathcal{C}^N$ the operator $\mathbf{A\otimes C}$ is
 a contraction. In other words,
 \[ \mathbf{A}\in C_{1,N}\cap L(\mathcal{X})^N \Longleftrightarrow
 L_\mathbf{A}\in\mathcal{AS}_N(\mathcal{X}),\]
 that coincides with  in \cite[Corollary~4.3]{K3} (here $L_\mathbf{A}(z):=z\mathbf{A},\ z\in\mathbb{C}^N$).
\end{rem}
Let us also note that using \cite[Corollary~4.3]{K3} one can deduce (v) from (i) directly. Indeed, if $\widetilde{\mathbf{A}}\in
L(\widetilde{\mathcal{X}})^N$ is a unitary $\rho$-dilation of
$\mathbf{A}\in L(\mathcal{X})^N$ then for any
$\mathbf{C}\in\mathcal{C}^N$ by \cite[Corollary~4.3]{K3} the operator
$\widetilde{\mathbf{A}}\otimes\mathbf{C}$ is a contraction.
Therefore, due to \cite{SzN},
$\widetilde{\mathbf{A}}\otimes\mathbf{C}\in
L(\widetilde{\mathcal{X}}\otimes\mathcal{H}_\mathbf{C})$ has a
unitary dilation $U\in L(\mathcal{K}),\
\mathcal{K}\supset\widetilde{\mathcal{X}}\otimes\mathcal{H}_\mathbf{C}$.
Then for any $n\in\mathbb{N}$:
\begin{eqnarray*}
 (\mathbf{A}\otimes\mathbf{C})^n &=&\rho
P_{\mathcal{X}\otimes\mathcal{H}_\mathbf{C}}(\widetilde{\mathbf{A}}\otimes\mathbf{C})^n|\mathcal{X}\otimes\mathcal{H}_\mathbf{C} \\
&=&
\rho P_{\mathcal{X}\otimes\mathcal{H}_\mathbf{C}}(
P_{\widetilde{\mathcal{X}}\otimes\mathcal{H}_\mathbf{C}}U^n|\widetilde{\mathcal{X}}\otimes\mathcal{H}_\mathbf{C})|
\mathcal{X}\otimes\mathcal{H}_\mathbf{C}\\
&=& \rho P_{\mathcal{X}\otimes\mathcal{H}_\mathbf{C}}U^n|
\mathcal{X}\otimes\mathcal{H}_\mathbf{C},
\end{eqnarray*}
 i.e., $U$ is a unitary
$\rho$-dilation of the operator $\mathbf{A}\otimes\mathbf{C}$.
Thus, $\mathbf{A}\otimes\mathbf{C}\in C_{\rho}$.

Let us define the \emph{numerical radius of an $N$-tuple of
operators} $\mathbf{A}\in L(\mathcal{X})^N$ as
\begin{equation}\label{eq:num-rad}
    w^{(N)}(\mathbf{A}):=\sup_{\mathbf{C}\in\mathcal{C}^N}w(\mathbf{A}\otimes\mathbf{C}).
\end{equation}
For $N=1$, $w^{(1)}(A)=w(A)$.  Indeed,
\begin{eqnarray*}
w^{(1)}(A) &=& \sup_{\| C\|<1}w(A\otimes
C) \geq \sup_{0<\varepsilon<1}w(A\otimes(1-\varepsilon)I_{\mathcal{H}_C})=\sup_{0<\varepsilon<1}(1-\varepsilon)w(A)\\
&=& w(A);\\
w^{(1)}(A)
 &=& \sup_{\| C\|<1}w(A\otimes C) \leq  \sup_{\|
C\|<1}w(A)\| C\|=w(A).
\end{eqnarray*}
Here we used the properties $w(A\otimes I_\mathcal{H})=w(A)$ and
$w(A\otimes B) \leq w(A)\| B\|$ valid for any $A\in
L(\mathcal{X}),B\in L(\mathcal{H})$ (see, e.g., \cite{FH}).
\begin{prop}\label{prop:n-num-rad}
$\mathbf{A}\in C_{2,N}\Longleftrightarrow w^{(N)}(\mathbf{A})\leq
1$.
\end{prop}
\begin{proof}
By Theorem~\ref{thm:n-char}, $\mathbf{A}\in C_{2,N}$ if and only
if $\mathbf{A\otimes C}\in C_2=C_{2,1}$ for any
$\mathbf{C}\in\mathcal{C}^N$. This, in turn, means that
$w(\mathbf{A\otimes C})\leq 1$ for any
$\mathbf{C}\in\mathcal{C}^N$ (by Berger's result mentioned in
Section~\ref{sec:prelim}), i.e, $w^{(N)}(\mathbf{A})\leq 1$.
\end{proof}
\begin{thm}\label{thm:pr}
If $\mathbf{A}\in C_{\rho,N}\cap L(\mathcal{X})^N$ for a $\rho>0$,
then $L_\mathbf{A}\in\rho\mathcal{AS}_N(\mathcal{X})$. For any
$\rho>0$ such that $\rho\neq 1$, there exists an $\mathbf{A}\in
L(\mathcal{X})^N$ such that
$L_\mathbf{A}\in\rho\mathcal{AS}_N(\mathcal{X})$ and
$\mathbf{A}\notin C_{\rho,N}$.
\end{thm}
\begin{proof}
Let $\mathbf{A}\in C_{\rho,N}\cap L(\mathcal{X})^N$ for some
$\rho>0$, and $\mathbf{C}\in\mathcal{C}^N$. Then $\mathbf{A}$ has
a unitary $\rho$-dilation $\widetilde{\mathbf{A}}\in
L(\widetilde{\mathcal{X}})^N$, and
\begin{eqnarray*}
\|\mathbf{A\otimes C}\| &=& \left\|\sum_{k=1}^NA_k\otimes
C_k\right\|=\left\|\left.\rho(P_\mathcal{X}\otimes
I_{\mathcal{H}_\mathbf{C}})\left(\sum_{k=1}^N\widetilde{A}_k\otimes
C_k\right)\right|\mathcal{X}\otimes\mathcal{H}_\mathbf{C}\right\| \\
&\leq &\rho\left\|\sum_{k=1}^N\widetilde{A}_k\otimes
C_k\right\|=\rho\left\|\widetilde{\mathbf{A}}\otimes\mathbf{C}\right\|\leq\rho
\end{eqnarray*}
(here we used again Corollary~4.3 in \cite{K3}). Thus,
$L_\mathbf{A}\in\rho\mathcal{AS}_N(\mathcal{X})$.

Now, let $0<\rho\neq 1$, and $\mathbf{A}\in L(\mathcal{X})^N$ be
such that
$\frac{1}{\rho}L_\mathbf{A}(\zeta)=\frac{1}{\rho}\zeta\mathbf{A}$
is a unitary operator for each $\zeta\in\mathbb{T}^N$. Then, again
by Corollary~4.3 in \cite{K3},
$L_\mathbf{A}\in\rho\mathcal{AS}_N(\mathcal{X})$. Suppose there
exists a unitary $\rho$-dilation $\widetilde{\mathbf{A}}\in
L(\widetilde{\mathcal{X}})^N$ of $\mathbf{A}$. Then for any
$\zeta\in\mathbb{T}^N$, $L_\mathbf{A}(\zeta)=\zeta\mathbf{A}=\rho
P_\mathcal{X}(\zeta\widetilde{\mathbf{A}})|\mathcal{X}$. Hence,
for any $\zeta\in\mathbb{T}^N$ and $x\in\mathcal{X}$,
\[
\|\zeta\widetilde{\mathbf{A}}x\|=\|x\|=\|\frac{1}{\rho}\zeta\mathbf{A}x\|=\|
P_\mathcal{X}(\zeta\widetilde{\mathbf{A}})x\|,\] that is possible
only if $\zeta\widetilde{\mathbf{A}}x\in\mathcal{X}$ for all
$\zeta\in\mathbb{T}^N$ and $x\in\mathcal{X}$. Therefore, for
$n>1$,
\[\rho^n\|x\|=\|(\zeta\mathbf{A})^nx\|=\|\rho
P_\mathcal{X}(\zeta\widetilde{\mathbf{A}})^nx\|=\rho\|(\zeta\widetilde{\mathbf{A}})^nx\|=\rho\|x\|,\]
that is impossible for $x\neq 0$. Thus, $\mathbf{A}\notin
C_{\rho,N}$.
\end{proof}
\begin{note*} Compare Theorem~\ref{thm:pr} with Remark~\ref{rem:1-contr}.
\end{note*}
The same argument as in the proof of the first part of
Theorem~\ref{thm:pr} shows that, for $\mathbf{A}\in C_{\rho,N}$,
\begin{equation}\label{eq:n-pb}
    \|(\mathbf{A\otimes C})^n\|\leq\rho,\quad n\in\mathbb{N},\
    \mathbf{C}\in\mathcal{C}^N.
\end{equation}
\begin{note*} Compare (\ref{eq:n-pb}) with (\ref{eq:1-pb}).
\end{note*}
This uniform (in $\mathbf{C}\in\mathcal{C}^N$) power-boundedness
of an $N$-tuple of operators $\mathbf{A}$ is, in our setting, a
generalization of power-boundedness of a single operator. Let us
define the \emph{spectral radius of an $N$-tuple of operators}
$\mathbf{A}\in L(\mathcal{X})^N$ as
\begin{equation}\label{eq:n-sp-rad}
    \nu^{(N)}(\mathbf{A}):=\lim_{n\rightarrow +\infty}\left(\sup_{\mathbf{C}\in\mathcal{C}^N}\|(\mathbf{A\otimes
    C})^n\|\right)^{\frac{1}{n}}.
\end{equation}
\begin{note*} Compare (\ref{eq:n-sp-rad}) with (\ref{eq:1-sp-rad}).
\end{note*}
In other words,
$\nu^{(N)}(\mathbf{A})=\nu^{(N,\infty)}(L_\mathbf{A})$, where
$\nu^{(N,\infty)}(f)$ is the \emph{spectral radius of an element
$f$ of the Banach algebra $H_N^\infty(\mathcal{X})$} consisting of
holomorphic $L(\mathcal{X})$-valued functions $f$ on
$\mathbb{D}^N$ which satisfy
\[ \| f\|_{\infty,N}:=\sup_{\mathbf{C}\in\mathcal{C}^N}\|
f(\mathbf{C})\|<\infty\] (this algebra was introduced in
\cite{Ag}). Here $f(\mathbf{C})$ is defined in the same manner as
$k(\mathbf{C,C})$ in (\ref{eq:n-calc}), i.e., for
\begin{eqnarray*}
f(z) &=& \sum_{t\in\mathbb{Z}^N_+}\hat{f}_tz^t,\quad
z\in\mathbb{D}^N,\\
f(\mathbf{C}) &:=& \sum_{t\in\mathbb{Z}^N_+}\hat{f}_t\otimes
\mathbf{C}^t,\quad \mathbf{C}\in\mathcal{C}^N,
\end{eqnarray*}
where the latter series converges in the
$L(\mathcal{X}\otimes\mathcal{H}_\mathbf{C})$-norm. For $N=1$,
$\nu^{(1)}(A)=\nu(A)$. Indeed,
\begin{eqnarray*}
\nu^{(1)}(A) &=& \lim_{n\rightarrow
+\infty}\left(\sup_{\|C\|<1}\|(A\otimes
    C)^n\|\right)^{\frac{1}{n}}=\lim_{n\rightarrow +\infty}\left(\sup_{\|C\|<1}\|A^n\otimes
    C^n\|\right)^{\frac{1}{n}}\\
    &=& \lim_{n\rightarrow +\infty}\left(\|A^n\|\sup_{\|C\|<1}
    \|C^n\|\right)^{\frac{1}{n}}=\lim_{n\rightarrow
    +\infty}\|A^n\|^{\frac{1}{n}}=\nu(A).
\end{eqnarray*}
\begin{rem}\label{rem:n-sp-rad}
For any $\mathbf{A}\in C_{\rho,N}$, by virtue of (\ref{eq:n-pb}),
$\nu^{(N)}(\mathbf{A})\leq 1$.
\end{rem}
\begin{thm}\label{thm:n-increase}
For a fixed Hilbert space $\mathcal{X}$ and any $N\geq 1$ the
class $C_{\rho,N}$ increases as a function of $\rho$:
\[ C_{\rho,N}\subset C_{\rho',N}\ \mbox{for}\ \rho<\rho'.\]
 Moreover, for
$\dim\mathcal{X}\geq 2$, $C_{\rho,N}$ increases strictly:
\[ C_{\rho,N}\neq C_{\rho',N}\ \mbox{for}\ \rho\neq\rho'.\]
For $\dim\mathcal{X}=1$ the classes $C_{\rho,N}$ coincide for all
$\rho\geq 1$, and strictly increase for $0<\rho<1$.
\end{thm}
\begin{proof}
For $N=1$ this theorem is true (see Section~\ref{sec:prelim}). For
$N>1$ it follows from the equivalence (i)$\Leftrightarrow$(v) in
Theorem~\ref{thm:n-char}.
\end{proof}
\begin{thm}\label{thm:n-vN}
For any $\mathbf{A}\in C_{\rho,N},\ \mathbf{C}\in\mathcal{C}^N$,
and a polynomial $p$ of one variable, \[\|
p(\mathbf{A}\otimes\mathbf{C})\|\leq\max_{|z|\leq 1}|\rho
p(z)+(1-\rho)p(0)|.\]
\end{thm}
\begin{proof}
This result follows from the generalized von Neumann inequality
(\ref{eq:vn}) and the equivalence (i)$\Leftrightarrow$(v) in
Theorem~\ref{thm:n-char}.
\end{proof}
Let us remark that results of this section on $N$-tuples of
operators from the classes $C_{\rho,N}$ can be extended to
elements of $H_N^\infty(\mathcal{X})$, though the notion of
unitary $\rho$-dilation no longer makes sense for this case. Define
$C_{\rho,N}^{(\infty)}$ as a class of functions $f\in
H_N^\infty(\mathcal{X})$ such that $f(\mathbf{C})\in
C_\rho=C_{\rho,1}$ for any $\mathbf{C}\in\mathcal{C}^N$. Then, in
particular, Theorem~\ref{thm:n-char} implies that $\mathbf{A}\in
C_{\rho,N}$ if and only if $L_\mathbf{A}\in
C_{\rho,N}^{(\infty)}$. The following analogue of
Theorem~\ref{thm:n-char} is easily obtained.
\begin{thm}\label{thm:n-char-bdd}
Let $f\in H_N^\infty(\mathcal{X})$ and $\rho>0$. The following
conditions are equivalent:
\begin{description}
    \item[(i)] $f\in C_{\rho,N}^{(\infty)}$;
    \item[(ii)] the function $k_{\rho,N}^f(z,w):=\rho
    I_\mathcal{X}-(\rho-1)(f(z)+(f(w)^*)+(\rho-2)f(w)^*f(z)$
is an Agler kernel on $\mathbb{D}^N\times\mathbb{D}^N$;
    \item[(iii)] the function
   $\psi_{\rho,N}^f(z):=(1-\frac{2}{\rho})I_\mathcal{X}+\frac{2}{\rho}(I_\mathcal{X}-f(z))^{-1}$
    belongs to $\mathcal{AH}_N(\mathcal{X})$;
    \item[(iv)] the function
    $\varphi_{\rho,N}^f(z):=f(z)((\varrho-1)f(z)-\rho I_\mathcal{X})^{-1}$
    belongs to $\mathcal{AS}_N(\mathcal{X})$.
\end{description}
\end{thm}
Clearly, $H_N^\infty(\mathcal{X})\cap
C_{1,N}^{(\infty)}=\mathcal{AS}_N(\mathcal{X})$. Set
\begin{equation}\label{eq:num-rad-bdd}
w^{(N,\infty)}(f):=\sup_{\mathbf{C}\in\mathcal{C}^N}w(f(\mathbf{C})).
\end{equation}
\begin{note*} Compare (\ref{eq:num-rad-bdd}) with (\ref{eq:num-rad}).
\end{note*}
\begin{rem}\label{rem:n-ext}
Proposition~\ref{prop:n-num-rad} extends directly  to the class
$C_{2,N}^{(\infty)}$, with $f\in H_N^\infty(\mathcal{X})$ in the
place of $\mathbf{A}\in L(\mathcal{X})^N$, and $w^{(N,\infty)}(f)$
in the place of $w^{(N)}(\mathbf{A})$. Remark~\ref{rem:n-sp-rad} extends directly 
to $f\in H_N^\infty(\mathcal{X})$ in the
place of $\mathbf{A}\in L(\mathcal{X})^N$, and
$\nu^{(N,\infty)}(f)$ in the place of $\nu^{(N)}(\mathbf{A})$.
Also, Theorems~\ref{thm:n-increase} and \ref{thm:n-vN} extend to the classes $C_{\rho,N}^{(\infty)}$.
\end{rem}

\section{Multivariable operator and operator-function radii}\label{sec:n-rad}
In this section we extend the notion of operator radii $w_\rho,\
0<\rho\leq\infty$, to the multivariable case, i.e., to $N$-tuples
of bounded linear operators and to elements of the Banach algebra
$H_N^\infty(\mathcal{X})$. Let $0<\rho<\infty$ and $f\in
H_N^\infty(\mathcal{X})$. Set
\[ w_{\rho,N}^{(\infty)}(f):=\inf\{ u>0:\frac{1}{u}f\in
C_{\rho,N}^{(\infty)}\},\] and for $\mathbf{A}\in
L(\mathcal{X})^N$, define
\[ w_{\rho,N}(\mathbf{A}):=w_{\rho,N}^{(\infty)}(L_\mathbf{A}).\]
Due to our remark preceding to Theorem~\ref{thm:n-char-bdd},
\begin{equation}\label{eq:n-rad}
w_{\rho,N}(\mathbf{A})=\inf\{ u>0:\frac{1}{u}\mathbf{A}\in
C_{\rho,N}\}.
\end{equation}
\begin{note*} Compare (\ref{eq:n-rad}) with (\ref{eq:1-rad}).
\end{note*}
Clearly, for $N=1$ and $A\in L(\mathcal{X})$,
$w_{\rho,1}(A)=w_\rho(A)$.
\begin{lem}\label{lem:w}
For $f\in H_N^\infty(\mathcal{X}),\ \mathbf{A}\in
L(\mathcal{X})^N$,
\begin{eqnarray}
w_{\rho,N}^{(\infty)}(f) &=&
\sup_{\mathbf{C}\in\mathcal{C}^N}w_\rho(f(\mathbf{C})),
\label{eq:w-f}\\
w_{\rho,N}(\mathbf{A}) &=&
\sup_{\mathbf{C}\in\mathcal{C}^N}w_\rho(\mathbf{A\otimes C}).
\label{eq:w-a}
\end{eqnarray}
\end{lem}
\begin{proof}
Let  $f\in H_N^\infty(\mathcal{X})$. Then for $u>0$,
$\frac{1}{u}f\in C_{\rho,N}^{(\infty)}$ if and only if for any
$\mathbf{C}\in\mathcal{C}^N$ one has $\frac{1}{u}f(\mathbf{C})\in
C_\rho$. Therefore,
\begin{eqnarray*}
w_{\rho,N}^{(\infty)}(f) &=& \inf\{ u>0:\frac{1}{u}f\in
C_{\rho,N}^{(\infty)}\}=\inf\{
u>0:\forall\mathbf{C}\in\mathcal{C}^N,\
\frac{1}{u}f(\mathbf{C})\in
C_\rho\}\\
&=& \sup_{\mathbf{C}\in\mathcal{C}^N}\inf\{
u>0:\frac{1}{u}f(\mathbf{C})\in C_\rho\}=
\sup_{\mathbf{C}\in\mathcal{C}^N}w_\rho(f(\mathbf{C})),
\end{eqnarray*}
i.e., (\ref{eq:w-f}) is true. Now, (\ref{eq:w-a}) follows from
(\ref{eq:w-f}) and the definition of $w_{\rho,N}(\mathbf{A})$.
\end{proof}
\begin{thm}\label{thm:n-rad}
{\bf 1.} All properties $(${\rm i}$)$--$(${\rm xx}$)$ listed in
Theorem~\ref{thm:1-rad} are satisfied for
$w_{\rho,N}^{(\infty)}(\cdot)$ in the place of $w_\rho(\cdot);\
f,g\in H_N^\infty(\mathcal{X})$ in the place of $A,B\in
L(\mathcal{X})$; $w^{(N,\infty)}(\cdot)$ in the place of
$w(\cdot)$; and $\nu^{(N,\infty)}(\cdot)$ in the place of
$\nu(\cdot)$.

 {\bf 2.} Properties $(${\rm i}$)$--$(${\rm
xvii}$)$ listed in Theorem~\ref{thm:1-rad} are satisfied for
$w_{\rho,N}(\cdot)$ in the place of $w_\rho(\cdot);\ \mathbf{A}\in
L(\mathcal{X})^N$ in the place of $A\in L(\mathcal{X})$;
$w^{(N)}(\cdot)$ in the place of $w(\cdot)$; and
$\nu^{(N)}(\cdot)$ in the place of $\nu(\cdot)$.
\end{thm}
\begin{proof}
{\bf 1.} Let $f\in H_N^\infty(\mathcal{X})$. Then $\|
f\|_{\infty,N}=\sup_{\mathbf{C}\in\mathcal{C}^N}\|
f(\mathbf{C})\|<\infty$.  By properties (vi) and (x) in
Theorem~\ref{thm:1-rad}, and Lemma~\ref{lem:w}, if $0<\rho\leq 1$
then
\begin{eqnarray*}
w_{\rho,N}^{(\infty)}(f) &=&
\sup_{\mathbf{C}\in\mathcal{C}^N}w_\rho(f(\mathbf{C}))\leq\left(\frac{2}{\rho}-1\right)
\sup_{\mathbf{C}\in\mathcal{C}^N}w_1(f(\mathbf{C})) \\
&=&
\left(\frac{2}{\rho}-1\right)\sup_{\mathbf{C}\in\mathcal{C}^N}\|
f(\mathbf{C})\| 
<\infty,
\end{eqnarray*}
 and if $\rho>1$ then
\[
w_{\rho,N}^{(\infty)}(f)=\sup_{\mathbf{C}\in\mathcal{C}^N}w_\rho(f(\mathbf{C}))\leq
\sup_{\mathbf{C}\in\mathcal{C}^N}w_1(f(\mathbf{C}))=
\sup_{\mathbf{C}\in\mathcal{C}^N}\|f(\mathbf{C})\|<\infty.\] Thus,
property (i) is fulfilled.

Properties (ii)--(vii), (ix)--(xi), (xiii)--(xv), and (xvii)--(xx)
easily follow from the properties in Theorem~\ref{thm:1-rad} with
the same numbers, and Lemma~\ref{lem:w}.

The proof of property (viii) is an adaptation of the proof of
Theorem~5.1 in \cite{H1} to our case. First of all, let us remark
that property (iv) implies that if $u>w_{\rho,N}^{(\infty)}(f)$
then $\frac{1}{u}f\in C_{\rho,N}^{(\infty)}$, and for any
$\mathbf{C}\in\mathcal{C}^N$ one has $\frac{1}{u}f(\mathbf{C})\in
C_\rho$. In particular,
\[
\sup_{\mathbf{C}\in\mathcal{C}^N}\left\|\left(\frac{f(\mathbf{C})}{u}\right)^n\right\|\leq\rho,\quad
n\in\mathbb{N}.\] Therefore,
$\nu^{(N,\infty)}\left(\frac{f}{u}\right)\leq 1$, i.e.,
$\nu^{(N,\infty)}(f)\leq u$. Thus, for any $\rho>0$,
$\nu^{(N,\infty)}(f)\leq w_{\rho,N}^{(\infty)}(f)$, moreover,
\[\nu^{(N,\infty)}(f)\leq \lim_{\rho\rightarrow
+\infty}w_{\rho,N}^{(\infty)}(f)\] (note, that due to property
(x), $w_{\rho,N}^{(\infty)}(f)$ is a non-increasing and bounded
from below function of $\rho$, hence it has a limit as
$\rho\rightarrow +\infty$).

For the proof of the opposite inequality, let us first show that
if $\nu^{(N,\infty)}(g)<1$ for some $g\in H_N^\infty(\mathcal{X})$
then beginning with some $\rho_0>0$ (i.e., for all $\rho\geq
\rho_0$) one has $g\in C_{\rho,N}^{(\infty)}$. Indeed, in  this 
 case there exists an $s>1$ such that
$\nu^{(N,\infty)}(sg)<1$. Then there exists a $B>0$ such that
\[ s^n\sup_{\mathbf{C}\in\mathcal{C}^N}\| g(\mathbf{C})^n\|\leq
B,\quad n\in\mathbb{N}.\] Hence, for any
$\mathbf{C}\in\mathcal{C}^N$,
\begin{eqnarray*}
\Real\psi_{\rho,N}^g(\mathbf{C}) &=&
\left(1-\frac{2}{\rho}\right)I_{\mathcal{X}\otimes\mathcal{H}_\mathbf{C}}+
\frac{2}{\rho}\Real(I_{\mathcal{X}\otimes\mathcal{H}_\mathbf{C}}-g(\mathbf{C}))^{-1}
\\ &=& I_{\mathcal{X}\otimes\mathcal{H}_\mathbf{C}}+
\frac{2}{\rho}\Real\sum_{n=1}^\infty g(\mathbf{C})^n \succeq
\left(1-\frac{2}{\rho}\sum_{n=1}^\infty \|
g(\mathbf{C})^n\|\right)I_{\mathcal{X}\otimes\mathcal{H}_\mathbf{C}}
\\
&\succeq & \left(1-\frac{2}{\rho}\sum_{n=1}^\infty
\frac{B}{s^n}\right)I_{\mathcal{X}\otimes\mathcal{H}_\mathbf{C}}=\left(1-\frac{2B}{\rho(s-1)}\right)I_{\mathcal{X}\otimes\mathcal{H}_\mathbf{C}}\succeq
0
\end{eqnarray*}
as soon as $\rho\geq\frac{2B}{s-1}$. Thus, by
Theorem~\ref{thm:n-char-bdd}, $g\in C_{\rho,N}^{(\infty)}$ for any
$\rho\geq\frac{2B}{s-1}$.

Now, if $\nu^{(N,\infty)}(f)=0$ then for any $k\in\mathbb{N}$,
$\nu^{(N,\infty)}(kf)=0$. Hence, for $\rho\geq\rho_0$ we have
$kf\in  C_{\rho,N}^{(\infty)}$, and by property (iv),
$w_{\rho,N}^{(\infty)}(kf)\leq 1$. Thus, $$\lim_{\rho\rightarrow
+\infty}w_{\rho,N}^{(\infty)}(f)\leq\frac{1}{k}$$ for any
$k\in\mathbb{N}$, and $$\lim_{\rho\rightarrow
+\infty}w_{\rho,N}^{(\infty)}(f)=0=\nu^{(N,\infty)}(f),$$ as
required.

If $\nu^{(N,\infty)}(f)>0$ then for any $\varepsilon>0$,
\[
\nu^{(N,\infty)}\left(\frac{f}{(1+\varepsilon)\nu^{(N,\infty)}(f)}\right)=\frac{1}{1+\varepsilon}<1.\]
Then for $\rho\geq\rho_0$,
\[
w_{\rho,N}^{(\infty)}\left(\frac{f}{(1+\varepsilon)\nu^{(N,\infty)}(f)}\right)\leq
1,\] hence
$w_{\rho,N}^{(\infty)}(f)\leq(1+\varepsilon)\nu^{(N,\infty)}(f)$.
Passing to the limit as $\rho\rightarrow+\infty$, and then as
$\varepsilon\downarrow 0$, we get
\[ \lim_{\rho\rightarrow
+\infty}w_{\rho,N}^{(\infty)}(f)\leq\nu^{(N,\infty)}(f),\] as
required. Thus, property (viii) is proved.

For the proof of property (xii), it is enough to suppose, by
virtue of positive homogeneity of $w_{\rho,N}^{(\infty)}(\cdot)$
and $\nu^{(N,\infty)}(\cdot)$, that for $f\in
H_N^\infty(\mathcal{X})$ one has $w_{\rho_0,N}^{(\infty)}(f)=1,\
\nu^{(N,\infty)}(f)<1$, and prove that for any $\rho>\rho_0$,
$w_{\rho,N}^{(\infty)}(f)<1$. By Theorem~\ref{thm:n-char-bdd} and
property (iv) in the present theorem,
\[
\frac{\rho_0}{2}\psi_{\rho_0,N}^f(z)=\left(\frac{\rho_0}{2}-1\right)I_\mathcal{X}+(I_\mathcal{X}-f(z))^{-1}\in\mathcal{AH}_N(\mathcal{X}),\]
i.e., for any $\mathbf{C}\in\mathcal{C}^N$,
\[
\Real\left[\frac{\rho_0}{2}\psi_{\rho_0,N}^f(\mathbf{C})\right]=\Real\left[\left(\frac{\rho_0}{2}-1\right)I_{\mathcal{X}\otimes\mathcal{H}_\mathbf{C}}
+(I_{\mathcal{X}\otimes\mathcal{H}_\mathbf{C}}-f(\mathbf{C}))^{-1}\right]\succeq
0,\] and for any $\rho>\rho_0$,
\[ \Real\left[\frac{\rho}{2}\psi_{\rho,N}^f(\mathbf{C})\right]=\Real\left[\left(\frac{\rho}{2}-1\right)I_{\mathcal{X}\otimes\mathcal{H}_\mathbf{C}}
+(I_{\mathcal{X}\otimes\mathcal{H}_\mathbf{C}}-f(\mathbf{C}))^{-1}\right]\succeq\frac{\rho-\rho_0}{2}I_{\mathcal{X}\otimes\mathcal{H}_\mathbf{C}}.\]
Since the resolvent $R_f(\lambda):=(\lambda
I_{\mathcal{X}}-f)^{-1}$ is
continuous in the $H_N^\infty(\mathcal{X})$-norm on the resolvent
set of $f$, and $\nu^{(N,\infty)}(f)<1$, for $\varepsilon>0$ small
enough, one has $\nu^{(N,\infty)}((1+\varepsilon)f)<1$, and
\[ \Real\left[\frac{\rho}{2}\psi_{\rho,N}^{(1+\varepsilon)f}(\mathbf{C})\right]=\Real\left[\left(\frac{\rho}{2}-1\right)
I_{\mathcal{X}\otimes\mathcal{H}_\mathbf{C}}
+(I_{\mathcal{X}\otimes\mathcal{H}_\mathbf{C}}-(1+\varepsilon)f(\mathbf{C}))^{-1}\right]\succeq
0\] for any $\mathbf{C}\in\mathcal{C}^N$, i.e.,
$\frac{\rho}{2}\psi_{\rho,N}^{(1+\varepsilon)f}\in\mathcal{AH}_N(\mathcal{X})$,
and
$\psi_{\rho,N}^{(1+\varepsilon)f}\in\mathcal{AH}_N(\mathcal{X})$.
Hence, by Theorem~\ref{thm:n-char-bdd}, $(1+\varepsilon)f\in
C_{\rho,N}^{(\infty)}$ which means, by property (iv), that
$w_{\rho,N}^{(\infty)}((1+\varepsilon)f)\leq 1$. Thus,
$w_{\rho,N}^{(\infty)}(f)\leq\frac{1}{1+\varepsilon}<1$, as
required.

The first part of property (xvi) in this theorem follows from
property (xvi) in Theorem~\ref{thm:1-rad}, and Lemma~\ref{lem:w}.
For the proof of the second part of (xvi), we use properties (xv)
and (xvi) from Theorem~\ref{thm:1-rad}, property (xv) in the
present theorem, and Lemma~\ref{lem:w}:
\begin{eqnarray*}
\lim_{\rho\downarrow 0}\frac{\rho}{2}w_{\rho,N}^{(\infty)}(f) &=&
\sup_{0<\rho<1}\left\{\frac{\rho}{2}w_{\rho,N}^{(\infty)}(f)\right\}=\sup_{0<\rho<1}\sup_{\mathbf{C}\in\mathcal{C}^N}
\left\{\frac{\rho}{2}w_\rho(f(\mathbf{C}))\right\}
\\
&=&
\sup_{\mathbf{C}\in\mathcal{C}^N}\sup_{0<\rho<1}\left\{\frac{\rho}{2}w_\rho(f(\mathbf{C}))\right\}=
\sup_{\mathbf{C}\in\mathcal{C}^N}\left\{\lim_{\rho\downarrow
0}\frac{\rho}{2}w_\rho(f(\mathbf{C}))\right\} \\
&=&
\sup_{\mathbf{C}\in\mathcal{C}^N}w_2(f(\mathbf{C}))=w_{2,N}^{(\infty)}(f).
\end{eqnarray*}
The proof of property (xvi), as well as part~\textbf{1} of this
theorem, is complete.

Part~{\bf 2} follows from part~{\bf 1}.
\end{proof}
Denote by $C_{\infty,N}^{(\infty)}$ (resp., $C_{\infty,N}$) the
class of $\mathcal{C}^N$-bounded holomorphic operator valued
functions on $\mathbb{D}^N$ (resp., $N$-tuples of bounded linear
operators on a common Hilbert space) with spectral radius at most
one.
\begin{thm}\label{thm:n-clos}
Let $\mathcal{X}$ be a Hilbert space. Then
\begin{eqnarray}
C_{\infty,N}^{(\infty)}\cap H_N^\infty(\mathcal{X})  &=&
\clos\left\{\bigcup_{0<\rho<\infty}(C_{\rho,N}^{(\infty)}\cap
H_N^\infty(\mathcal{X}))\right\}; \label{eq:sp-r-f} \\
C_{\infty,N}\cap L(\mathcal{X})^N  &=&
\clos\left\{\bigcup_{0<\rho<\infty}(C_{\rho,N}\cap
L(\mathcal{X})^N)\right\}. \label{eq:sp-r-t}
\end{eqnarray}
\end{thm}
\begin{note*}
Compare (\ref{eq:sp-r-f}) and (\ref{eq:sp-r-t}) with
(\ref{eq:1-clos}).
\end{note*}
\begin{proof}[Proof of Theorem~\ref{thm:n-clos}]
The inclusion ``$\supset$" in (\ref{eq:sp-r-f}) and
(\ref{eq:sp-r-t}) follows from Remarks \ref{rem:n-sp-rad} and
\ref{rem:n-ext}, and the fact that the set of
$\mathcal{C}^N$-bounded holomorphic operator-valued functions on
$\mathbb{D}^N$ (resp., $N$-tuples of bounded operators) with
spectral radius at most one is closed in $H_N^\infty(\mathcal{X})$
(resp., $L(\mathcal{X})^N$).

To show the inclusion ``$\subset$" in (\ref{eq:sp-r-f}), observe
that for $f\in C_{\infty,N}^{(\infty)}\cap
H_N^\infty(\mathcal{X})$ and $0<r<1$, $\nu^{(N,\infty)}(rf)\leq
r<1$. By property (viii) from Theorem~\ref{thm:n-rad}, for
$\rho_0>0$ big enough, $w_{\rho_0,N}^{(\infty)}(rf)<1$, and by
property (iv) from the same theorem,
\[ rf\in C_{\rho_0,N}^{(\infty)}\cap
H_N^\infty(\mathcal{X})\subset\clos\left\{\bigcup_{0<\rho<\infty}(C_{\rho,N}^{(\infty)}\cap
H_N^\infty(\mathcal{X}))\right\}.\] Passing to the limit as
$r\uparrow 1$, we get
\[ f\in\clos\left\{\bigcup_{0<\rho<\infty}(C_{\rho,N}\cap
H_N^\infty(\mathcal{X}))\right\},\] and the inclusion ``$\subset$"
in (\ref{eq:sp-r-f}) follows. Analogously for the inclusion
``$\subset$" in (\ref{eq:sp-r-t}).
\end{proof}
In view of property (iv) in Theorem~\ref{thm:n-rad}, let us call
the elements of the class $C_{\rho,N}$ \emph{($N$-variable)
$\rho$-contractions}.

\section{On similarity of $\rho$-contractions to $1$-contractions in several variables}\label{sec:similar}
An $N$-tuple of operators $\mathbf{A}=(A_1,\ldots,A_N)\in
L(\mathcal{X})^N$ is said to be \emph{simultaneously similar} to
an $N$-tuple of operators $\mathbf{B}=(B_1,\ldots,B_N)\in
L(\mathcal{Y})^N$ if there exists a boundedly invertible operator
$S\in L(\mathcal{X,Y})$ such that
\begin{equation}\label{eq:n-similar}
    A_k=S^{-1}B_kS,\quad k=1,\ldots,N,
\end{equation}
or equivalently,
\begin{equation}\label{eq:n-pencil-similar}
    z\mathbf{A}=S^{-1}(z\mathbf{B})S,\quad z\in\mathbb{C}^N.
\end{equation}
\begin{note*}
Compare (\ref{eq:n-similar}) and (\ref{eq:n-pencil-similar}) with
(\ref{eq:1-similar}).
\end{note*}
\begin{thm}\label{thm:non-similar}
For any $\rho>1$ and $N>1$, there exists an
$\mathbf{A}=(A_1,\ldots,A_N)\in C_{\rho,N}$ which is not
simultaneously similar to any $\mathbf{T}=(T_1,\ldots,T_N)\in
C_{1,N}$.
\end{thm}
\begin{proof}
Let $N=2$, and for any $\varepsilon\geq 0$ set
$\mathbf{A}^{(\varepsilon)}=(A_1^{(\varepsilon)},A_2^{(\varepsilon)})\in
L(\mathbb{C}^3)^2$, where
\[ A_1^{(\varepsilon)}:=\left[\begin{array}{lll}
  0 & \frac{1+\varepsilon}{\sqrt{2}} & 0 \\
  0 & 0 & 0 \\
  -\frac{1+\varepsilon}{\sqrt{2}} & 0 & 0
\end{array}\right],\quad  A_2^{(\varepsilon)}:=\left[\begin{array}{ccc}
  0 & 0 & \frac{1+\varepsilon}{\sqrt{2}} \\
  \frac{1+\varepsilon}{\sqrt{2}} & 0 & 0 \\
  0 & 0 & 0
\end{array}\right].  \]
Then for any $\varepsilon\geq 0$ and $z\in\mathbb{C}^2$,
\begin{eqnarray*}
z\mathbf{A}^{(\varepsilon)} &=& \left[\begin{array}{ccc}
  0 & \frac{1+\varepsilon}{\sqrt{2}}z_1 & \frac{1+\varepsilon}{\sqrt{2}}z_2 \\
  \frac{1+\varepsilon}{\sqrt{2}}z_2 & 0 & 0 \\
  -\frac{1+\varepsilon}{\sqrt{2}}z_1 & 0 & 0
\end{array}\right], \\ (z\mathbf{A}^{(\varepsilon)})^2 &=& \left[\begin{array}{ccc}
  0 & 0 & 0 \\
  0 & \frac{(1+\varepsilon)^2}{2}z_1z_2 & \frac{(1+\varepsilon)^2}{2}z_2^2 \\
  0 & -\frac{(1+\varepsilon)^2}{2}z_1^2 & -\frac{(1+\varepsilon)^2}{2}z_1z_2
\end{array}\right],\\
(z\mathbf{A}^{(\varepsilon)})^3 &=&
(z\mathbf{A}^{(\varepsilon)})^4=\ldots=0,
\end{eqnarray*}
i.e., $z\mathbf{A}^{(\varepsilon)}$ is a nilpotent operator of
degree 3. Hence, for any $\rho>1$ and $z\in\mathbb{D}^2$,
\begin{eqnarray*}
\|\varphi_{\rho,2}^{\mathbf{A}^{(0)}}(z)\| &=& \left\|
z\mathbf{A}^{(0)}((\rho-1)z\mathbf{A}^{(0)}-\rho
I)^{-1}\right\|=\left\|\frac{z\mathbf{A}^{(0)}}{\rho}\left(I-\frac{\rho-1}{\rho}z\mathbf{A}^{(0)}\right)^{-1}\right\|
\\
&=&
\left\|\frac{z\mathbf{A}^{(0)}}{\rho}+(\rho-1)\left(\frac{z\mathbf{A}^{(0)}}{\rho}\right)^2\right\|\leq\frac{1}{\rho}\|
z\mathbf{A}^{(0)}\|+\frac{\rho-1}{\rho^2}\|
(z\mathbf{A}^{(0)})^2\| \\
&=& \frac{1}{\rho}\left\|\left[\begin{array}{ccc}
  0 & \frac{z_1}{\sqrt{2}} & \frac{z_2}{\sqrt{2}} \\
  \frac{z_2}{\sqrt{2}} & 0 & 0 \\
  -\frac{z_1}{\sqrt{2}} & 0 & 0
\end{array}\right]\right\| +\frac{\rho-1}{\rho^2}\left\|\left[\begin{array}{c}
  0 \\
  \frac{z_2}{\sqrt{2}} \\
 -\frac{z_1}{\sqrt{2}}
\end{array}\right]\left[\begin{array}{c}
  0 \\
  \frac{z_1}{\sqrt{2}} \\
 \frac{z_2}{\sqrt{2}}
\end{array}\right]^T\right\| \\
&\leq & \frac{1}{\rho}
+\frac{\rho-1}{\rho^2}=\frac{2\rho-1}{\rho^2}<1.
\end{eqnarray*}
Then, due to the von Neumann inequality in two variables
\cite{An}, one has
\[\|\varphi_{\rho,2}^{\mathbf{A}^{(0)}}(\mathbf{C})\|\leq\frac{2\rho-1}{\rho^2}<1,
\quad \mathbf{C}\in\mathcal{C}^2,\] i.e.,
$\varphi_{\rho,2}^{\mathbf{A}^{(0)}}\in\mathcal{AS}_2(\mathbb{C}^3)$.
Analogously, for $\varepsilon>0$ small enough (the choice of
$\varepsilon$ depends on $\rho$), one has
\[\sup_{\mathbf{C}\in\mathcal{C}^2}\|\varphi_{\rho,2}^{\mathbf{A}^{(\varepsilon)}}(\mathbf{C})\|=
\sup_{z\in\mathbb{D}^2}\|\varphi_{\rho,2}^{\mathbf{A}^{(\varepsilon)}}(\mathbf{z})\|\leq\frac{1+\varepsilon}{\rho}+
(\rho-1)\left(\frac{1+\varepsilon}{\rho}\right)^2<1,\] i.e.,
$\varphi_{\rho,2}^{\mathbf{A}^{(\varepsilon)}}\in\mathcal{AS}_2(\mathbb{C}^3)$,
and by Theorem~\ref{thm:n-char}, $\mathbf{A}^{(\varepsilon)}\in
C_{\rho,2}$.

Let us show now that for any $\varepsilon>0$ the pair
$\mathbf{A}^{(\varepsilon)}=(A_1^{(\varepsilon)},A_2^{(\varepsilon)})$
is not simultaneously similar to any pair $\mathbf{T}=(T_1,T_2)\in
C_{1,2}$. Observe that
\begin{eqnarray*}
(A_1^{(\varepsilon)}+A_2^{(\varepsilon)})(A_1^{(\varepsilon)}-A_2^{(\varepsilon)})
&=& \left[\begin{array}{ccc}
  0 & \frac{1+\varepsilon}{\sqrt{2}} & \frac{1+\varepsilon}{\sqrt{2}} \\
  \frac{1+\varepsilon}{\sqrt{2}} & 0 & 0 \\
  -\frac{1+\varepsilon}{\sqrt{2}} & 0 & 0
\end{array}\right]\left[\begin{array}{ccc}
  0 & \frac{1+\varepsilon}{\sqrt{2}} & -\frac{1+\varepsilon}{\sqrt{2}} \\
  -\frac{1+\varepsilon}{\sqrt{2}} & 0 & 0 \\
  -\frac{1+\varepsilon}{\sqrt{2}} & 0 & 0
\end{array}\right] \\
&=& \left[\begin{array}{ccc}
  -(1+\varepsilon)^2 & 0 & 0 \\
  0 & \frac{(1+\varepsilon)^2}{2} & -\frac{(1+\varepsilon)^2}{2} \\
  0 & -\frac{(1+\varepsilon)^2}{2} & \frac{(1+\varepsilon)^2}{2}
\end{array}\right].
\end{eqnarray*}
Then
\begin{equation}\label{eq:lim}
    \lim_{n\rightarrow +\infty}\|
[(A_1^{(\varepsilon)}+A_2^{(\varepsilon)})(A_1^{(\varepsilon)}-A_2^{(\varepsilon)})]^n\|=\infty.
\end{equation}
On the other hand, if
$\mathbf{A}^{(\varepsilon)}=(A_1^{(\varepsilon)},A_2^{(\varepsilon)})$
is simultaneously similar to some  $\mathbf{T}=(T_1,T_2)\in
C_{1,2}$ then for any $n\in\mathbb{N}$ one would have
\[\|
[(A_1^{(\varepsilon)}+A_2^{(\varepsilon)})(A_1^{(\varepsilon)}-A_2^{(\varepsilon)})]^n\|=\|
S^{-1}[(T_1+T_2)(T_1-T_2)]^nS\|\leq\| S\|\|S^{-1}\|<\infty,\]
since $\| T_1\pm T_2\|\leq 1$. We get a contradiction with
(\ref{eq:lim}).

Examples  of $N$-tuples of operators from $C_{\rho,N},\ \rho>1$,
which are not simultaneously similar to any $\mathbf{T}\in
C_{1,N}$ for the case $N>2$ can be obtained from the examples of
pairs $\mathbf{A}=(A_1^{(\varepsilon)},A_2^{(\varepsilon)})$
above, for sufficiently small $\varepsilon>0$, by setting zeros
for the rest of operators in these $N$-tuples:
$\widetilde{\mathbf{A}}^{(\varepsilon)}:=(A_1^{(\varepsilon)},A_2^{(\varepsilon)},0,\ldots,0)$.
\end{proof}
Let $\mathbf{A}=(A_1,\ldots,A_N)\in L(\mathcal{X})^N$. Then
$\widetilde{\mathbf{A}}=(\widetilde{A}_1,\ldots,\widetilde{A}_N)\in
L(\widetilde{\mathcal{X}})^N$ is called a \emph{uniform
$\rho$-dilation} of $\mathbf{A}$ if
$\widetilde{\mathcal{X}}\supset\mathcal{X}$ and
\begin{equation}\label{eq:uniform}
\forall n\in\mathbb{N},\ \forall i_1,\ldots,i_n\in\{
1,\ldots,N\},\quad A_{i_1}\cdots A_{i_n}=\rho
P_\mathcal{X}\widetilde{A}_{i_1}\cdots\widetilde{A}_{i_n}|\mathcal{X},
\end{equation}
or equivalently,
\begin{equation}\label{eq:uniform'}
\forall n\in\mathbb{N},\ \forall
z^{(1)},\ldots,z^{(n)}\in\mathbb{C}^N,\quad
z^{(1)}\mathbf{A}\cdots z^{(n)}\mathbf{A}=\rho
P_\mathcal{X}z^{(1)}\widetilde{\mathbf{A}}\cdots
z^{(n)}\widetilde{\mathbf{A}}|\mathcal{X}.
\end{equation}
\begin{note*}
Compare (\ref{eq:uniform}) and (\ref{eq:uniform'}) with
(\ref{eq:n-rhodil}) and (\ref{eq:n-rhodil'}).
\end{note*}
Clearly, a uniform $\rho$-dilation is a $\rho$-dilation. If
$\widetilde{\mathbf{A}}\in L(\widetilde{\mathcal{X}})^N$ is a
uniform $\rho$-dilation of $\mathbf{A}\in L(\mathcal{X})$, and for
any $\zeta\in\mathbb{T}^N$, $\zeta\widetilde{\mathbf{A}}$ is a
unitary operator, then $\widetilde{\mathbf{A}}$ is called a
\emph{uniform unitary $\rho$-dilation}  of $\mathbf{A}$. Denote by
$C_{\rho,N}^u$ the class of $N$-tuples of operators
$\mathbf{A}=(A_1,\ldots,A_N)$ on a common Hilbert space which
admit a uniform unitary $\rho$-dilation. Clearly,
$C_{\rho,N}^u\subset C_{\rho,N}$.
\begin{thm}\label{thm:similar}
Any $\mathbf{A}=(A_1,\ldots,A_N)\in C_{\rho,N}^u$ is
simultaneously similar to some $\mathbf{T}=(T_1,\ldots,T_N)\in
C_{1,N}^u$.
\end{thm}
\begin{proof}
Let $\mathbf{A}=(A_1,\ldots,A_N)\in C_{\rho,N}^u\cap
L(\mathcal{X})^N$, and $\mathbf{U}=(U_1,\ldots,U_N)\in
L(\widetilde{\mathcal{X}})^N$ be a uniform unitary $\rho$-dilation
of $\mathbf{A}$. Let $\mathcal{A}\subset
L(\widetilde{\mathcal{X}})$ be the minimal $C^*$-algebra which
contains the operators
$I_{\widetilde{\mathcal{X}}},U_1,\ldots,U_N$, and
$\mathcal{B}\subset L(\widetilde{\mathcal{X}})$ be the minimal
algebra over $\mathbb{C}$ which contains the operators
$U_1,\ldots,U_N$. Clearly, $\mathcal{B}\subset\mathcal{A}$. Let
$\varphi:\mathcal{B}\longrightarrow L(\mathcal{X})$ be a
homomorphism defined on the generators as
\[\varphi: U_k\longmapsto A_k,\ k=1,\ldots,N.\]
The algebra $\mathcal{B}$ consists of operators of the form
\[ p(\mathbf{U})=\sum_{1\leq k\leq m,\ i_1,\ldots,i_k\in\{ 1,\ldots,N\}
}\alpha_{i_1,\ldots,i_k}U_{i_1}\cdots U_{i_k},\] where
$\alpha_{i_1,\ldots,i_k}\in\mathbb{C}$ for all
$i_1,\ldots,i_k\in\{ 1,\ldots,N\}$. Then
\begin{eqnarray*}
\varphi(p(\mathbf{U})) &=&
\varphi(\sum\alpha_{i_1,\ldots,i_k}U_{i_1}\cdots
U_{i_k})=\sum\alpha_{i_1,\ldots,i_k}A_{i_1}\cdots
A_{i_k}\\
&=& p(\mathbf{A})=\rho P_\mathcal{X}p(\mathbf{U})|\mathcal{X}.
\end{eqnarray*}
Therefore, if $p(\mathbf{U})=0$ then $\varphi(p(\mathbf{U}))=0$,
and $\varphi$ is correctly defined. The homomorphism $\varphi$ is
\emph{completely bounded}, i.e.,
\[\|\varphi\|_{cb}:=\sup_{n\in\mathbb{N}}\|\id_n\otimes\varphi\|
<\infty,\] where $\id_n$ is the identical map of the matrix
algebra $\mathcal{M}_n(\mathbb{C})$ onto itself. Moreover,
$\|\varphi\|_{cb}\leq\rho$. Indeed, for any $n\in\mathbb{N}$ and a
polynomial $n\times n$ matrix of $N$ non-commuting variables,
\[ P(\mathbf{X})=[p_{ij}(\mathbf{X})]_{i,j=1}^n=\left[\sum_{1\leq k\leq m,\ i_1,\ldots,i_k\in\{ 1,\ldots,N\}
}\alpha_{i_1,\ldots,i_k}^{(ij)}X_{i_1}\cdots
X_{i_k}\right]_{i,j=1}^n,\]
\begin{eqnarray*}
\|(\id_n\otimes\varphi)(P(\mathbf{U}))\| &=&
\|(\id_n\otimes\varphi)\left([p_{ij}(\mathbf{U})]_{i,j=1}^n\right)\|
=\|[\varphi(p_{ij}(\mathbf{U}))]_{i,j=1}^n\| \\
&=& \|[p_{ij}(\mathbf{A})]_{i,j=1}^n\|=\|[\rho P_\mathcal{X}p_{ij}(\mathbf{U})|\mathcal{X}]_{i,j=1}^n\| \\
&=& \rho\|(I_{\mathbb{C}^n}\otimes
P_\mathcal{X})[p_{ij}(\mathbf{U})]_{i,j=1}^n|\mathbb{C}^n\otimes\mathcal{X}\|
\\
&\leq & \rho\| [p_{ij}(\mathbf{U})]_{i,j=1}^n\| = \rho\|
P(\mathbf{U})\|.
\end{eqnarray*}
Then, by Theorem~3.1 in \cite{Pa}, there exist a Hilbert space
$\mathcal{N}$, a \emph{completely contractive} homomorphism
$\gamma:\mathcal{B}\longrightarrow L(\mathcal{N})$ (i.e., such
that $\|\gamma\|_{cb}\leq 1$), and a boundedly invertible operator
$S\in L(\mathcal{X},\mathcal{N})$ such that
\[ \varphi(b)=S^{-1}\gamma(b)S,\quad b\in\mathcal{B}.\]
Moreover, as was shown in the proof of Theorem~3.1 in \cite{Pa},
$\gamma$ can be chosen in the form
\[ \gamma(b)=P_\mathcal{N}\pi(b)|\mathcal{N},\quad b\in\mathcal{B},\]
where $\pi:\mathcal{A}\longrightarrow L(\mathcal{K})$ is a
$*$-homomorphism, for some Hilbert space
$\mathcal{K}\supset\mathcal{N}$. In addition, it follows from
Theorem~2.7 and the proof of Theorem~2.8 in \cite{Pa} that one can
choose $\mathcal{K}=\mathcal{K}_1\oplus\mathcal{K}_1$, for some
Hilbert space $\mathcal{K}_1$, and
\[ \pi(a)=\pi_1(a)\oplus 0,\quad a\in\mathcal{A},\]
where $\pi_1:\mathcal{A}\longrightarrow L(\mathcal{K}_1)$ is a
unital $*$-homomorphism. Set \[ T_k:=\gamma(U_k)\in
L(\mathcal{N}),\quad k=1,\ldots,N.\] Then
\[ A_k=\varphi(U_k)=S^{-1}T_kS,\quad k=1,\ldots,N.\]
It remains to show that $\mathbf{T}=(T_1,\ldots,T_N)\in
C_{1,N}^u$. Set \[W_k:=\pi(U_k)\in L(\mathcal{K}),\quad
k=1,\ldots,N.\]  Since for any $n\in\mathbb{N}$ and
$i_1,\ldots,i_n\in\{ 1,\ldots,N\}$ one has
\begin{eqnarray*}
T_{i_1}\cdots T_{i_n} &=&
\gamma(U_{i_1}\cdots
U_{i_n})=P_\mathcal{N}\pi(U_{i_1}\cdots
U_{i_n})|\mathcal{N}\\
&=& P_\mathcal{N}W_{i_1}\cdots W_{i_n}|\mathcal{N},
\end{eqnarray*}
$\mathbf{W}=(W_1,\ldots,W_N)$ is a uniform $1$-dilation of
$\mathbf{T}=(T_1,\ldots,T_N)$, however, still not unitary.
Actually,
\[ W_k=\pi(U_k)=\pi_1(U_k)\oplus 0\ (=: W_k^{(1)}\oplus 0),\quad  k=1,\ldots,N.\]
Since $\pi_1$ is a unital $*$-homomorphism, and for any
$\zeta\in\mathbb{T}^N$,
\[
(\zeta\mathbf{U})^*\zeta\mathbf{U}=I_{\widetilde{\mathcal{X}}}=\zeta\mathbf{U}(\zeta\mathbf{U})^*,\]
one has, for any $\zeta\in\mathbb{T}^N$,
\[
(\zeta\mathbf{W}^{(1)})^*\zeta\mathbf{W}^{(1)}=I_{\mathcal{K}_1}=\zeta\mathbf{W}^{(1)}(\zeta\mathbf{W}^{(1)})^*,\]
where $\mathbf{W}^{(1)}=(W_1^{(1)},\ldots,W_N^{(1)})$. Set
\[ \widetilde{W}_k:=W_k^{(1)}\oplus\delta_{1k}V\in
L(\mathcal{K}_1\oplus\mathcal{R}),\quad k=1,\ldots,N,\] where
$\delta_{ij}$ is the Kronecker symbol, and $V$ is a unitary
dilation of the zero operator on $\mathcal{K}_1$, e.g., the
two-sided shift on the space
$\mathcal{R}:=l^2(\mathcal{K}_1)=\bigoplus_{-\infty}^{+\infty}\mathcal{K}_1$
(here we identify the space $\mathcal{K}_1$ with the subspace
$\ldots\oplus\{ 0\}\oplus\{ 0\}\oplus\mathcal{K}_1\oplus\{
0\}\oplus\{ 0\}\oplus\ldots$ in $\mathcal{R}$). Then, for any
$\zeta\in\mathbb{T}^N$,
\[
(\zeta\widetilde{\mathbf{W}})^*\zeta\widetilde{\mathbf{W}}=I_{\mathcal{K}_1\oplus\mathcal{R}}=
\zeta\widetilde{\mathbf{W}}(\zeta\widetilde{\mathbf{W}})^*,\] and
$\widetilde{\mathbf{W}}=(\widetilde{W}_1,\ldots,\widetilde{W}_N)\in
L(\mathcal{K}_1\oplus\mathcal{R})$ is a uniform unitary
$1$-dilation of $\mathbf{W}=(W_1,\ldots,W_N)$, and therefore, of
$\mathbf{T}=(T_1,\ldots,T_N)$. Thus, $\mathbf{T}\in C_{1,N}^u$, as
required.
\end{proof}
Theorem~\ref{thm:similar} is similar to the result of
G.~Popescu \cite{Po2} on simultaneous similarity of
$\rho$-contractions to  $1$-contractions in several variables,
however his notion of multivariable $\rho$-contractions is
different. Let us clarify the relation between these two results.
Denote by $C_{\rho,N}^P$ (we use here this notation instead of
just $C_\rho$, as in \cite{Po2}) the \emph{Popescu class} of all
$N$-tuples $\mathbf{A}=(A_1,\ldots,A_N)$ of bounded linear
operators on a common Hilbert space, say $\mathcal{X}$, which have
a \emph{uniform isometric $\rho$-dilation}, i.e., such an
$N$-tuple of operators $\mathbf{V}=(V_1,\ldots,V_N)\in
L(\widetilde{\mathcal{X}})^N,\
\widetilde{\mathcal{X}}\supset\mathcal{X}$, for which
\begin{description}
    \item[(1)] $V_k^*V_k=I_{\widetilde{\mathcal{X}}},\quad k=1,\ldots,N$;
\item[(2)] $V_k^*V_j=0,\quad k\neq j$;
\item[(3)] $\forall n\in\mathbb{N},\ \forall i_1,\ldots,i_n,\quad A_{i_1}\cdots A_{i_n}=\rho P_\mathcal{X}V_{i_1}\cdots
V_{i_n}|\mathcal{X}.$
\end{description}
Condition (2) in this definition can be replaced by
\begin{description}
\item[(2')] $\sum\limits_{k=1}^NV_kV_k^*\preceq
I_{\widetilde{\mathcal{X}}}$,
\end{description}
\noindent since (1)$\&$(2)$\Longleftrightarrow$(1)$\&$(2').
According to \cite{Po1}, the class $C_{1,N}^P$ coincides with the
class of $N$-tuples of operators $\mathbf{A}=(A_1,\ldots,A_N)$
($\in L(\mathcal{X})^N$, for some Hilbert space $\mathcal{X}$)
such that \[ \sum_{k=1}^NA_kA_k^*\preceq I_\mathcal{X}.\] By
Theorem~4.5 in \cite{Po2}, any $\mathbf{A}=(A_1,\ldots,A_N)\in
C_{\rho,N}^P,\ \rho>0$, is simultaneously similar to some
$\mathbf{T}=(T_1,\ldots,T_N)\in C_{1,N}^P$. This is a generalization of the
theorem of Sz.-Nagy and Foia\c{s} \cite{SzNF2} to several
variables. Theorem~\ref{thm:similar} of the present paper is a different generalization of the same result,
 since our classes
$C_{\rho,N}^u$ are different from Popescu's classes $C_{\rho,N}^P$
for $N>1$. More precisely, the following is true.
\begin{thm}\label{thm:pop}
For any $N>1$ and $\rho>0$, $C_{\rho,N}^u\varsubsetneq
C_{\rho,N}^P$.
\end{thm}
\begin{proof}
Let $\mathbf{A}\in C_{\rho,N}^u\cap L(\mathcal{X})^N,\ N>1$, and
$\mathbf{U}\in L(\widetilde{\mathcal{X}})^N$ be a uniform unitary
$\rho$-dilation of $\mathbf{A}$. Since for any
$\zeta\in\mathbb{T}^N$ the operator $\zeta\mathbf{U}$ is unitary,
it follows that
\[
\zeta\mathbf{U}(\zeta\mathbf{U})^*=I_{\widetilde{\mathcal{X}}},\quad\zeta\in\mathbb{T}^N,\]
which implies
\[\sum_{k=1}^NU_kU_k^*=I_{\widetilde{\mathcal{X}}}.\]
Thus, by \cite{Po1}, $\mathbf{U}\in C_{1,N}^P$. Let $\mathbf{V}\in
L(\hat{\mathcal{X}})^N$ be a uniform isometric $1$-dilation of
$\mathbf{U}$ in the sense of Popescu. Then for any
$n\in\mathbb{N}$ and $i_1,\ldots,i_n\in\{ 1,\ldots,N\}$,
\begin{eqnarray*}
A_{i_1}\cdots A_{i_n} &=& \rho
P_\mathcal{X}U_{i_1}\cdots U_{i_n}|\mathcal{X}=\rho 
P_\mathcal{X}(P_{\widetilde{\mathcal{X}}}V_{i_1}\cdots
V_{i_n}|\widetilde{\mathcal{X}})|\mathcal{X} \\
&=& \rho P_\mathcal{X}V_{i_1}\cdots
V_{i_n}|\mathcal{X},
\end{eqnarray*}
i.e., $\mathbf{V}$ is a uniform isometric $\rho$-dilation of
$\mathbf{A}$ in the sense of Popescu. Thus, $\mathbf{A}\in
C_{\rho,N}^P$. This proves the inclusion $C_{\rho,N}^u\subset
C_{\rho,N}^P$.

Let us prove that this inclusion is proper for any $N>1$ and
$\rho>0$. Firstly, consider the case $N=2$. Let $B\in
L(\mathcal{X}_0)$ be any operator of the class $C_\rho$ with $\|
B\|=\rho$. For example,
\[ B:=\left[\begin{array}{cc}
  0 & \rho \\
  0 & 0
\end{array}\right]\in L(\mathbb{C}^2)\]
satisfies $B^2=0$ and $\| B\|=\rho$, therefore by properties (iii)
and (xi) in Theorem~\ref{thm:1-rad}, $w_\rho(B)=1$, and by
property (iv) in the same theorem, $B\in C_\rho$. Set
$\mathcal{X}:=\mathcal{X}_0\oplus\mathcal{X}_0$,
\begin{equation}\label{eq:a's}
    A_1:=\left[\begin{array}{cc}
  B & 0 \\
  0 & 0
\end{array}\right]\in L(\mathcal{X}),\quad  A_2:=\left[\begin{array}{cc}
  0 & 0 \\
  B & 0
\end{array}\right]\in L(\mathcal{X}).
\end{equation}
Let $U\in L(\widetilde{\mathcal{X}}_0)$ be a unitary
$\rho$-dilation of $B$. Set
$\widetilde{\mathcal{X}}:=\widetilde{\mathcal{X}}_0\oplus\widetilde{\mathcal{X}}_0\oplus\ldots$,
and identify $\mathcal{X}=\mathcal{X}_0\oplus\mathcal{X}_0$ with
the subspace $\mathcal{X}_0\oplus\mathcal{X}_0\oplus\{ 0\}\oplus\{
0\}\oplus\ldots$ in $\widetilde{\mathcal{X}}$. Set
\[ V_1:=\left[\begin{array}{ccc}
  \fbox{$\begin{array}{c}
    U \\
    0
  \end{array}$} &  &  \\
   & \fbox{$\begin{array}{c}
    U \\
    0
  \end{array}$} &  \\
   &  & \ddots \\
\end{array}\right]\in L(\widetilde{\mathcal{X}}),\quad  V_2:=\left[\begin{array}{ccc}
  \fbox{$\begin{array}{c}
    0 \\
    U
  \end{array}$} &  &  \\
   & \fbox{$\begin{array}{c}
    0 \\
    U
  \end{array}$} &  \\
   &  & \ddots \\
\end{array}\right]\in L(\widetilde{\mathcal{X}}),\]
i.e., the operators $V_1$ and $V_2$ are introduced here as
infinite block-diagonal matrices with equal operator blocks
$\left[\begin{array}{c}
    U \\
    0
  \end{array}\right]\in
L(\widetilde{\mathcal{X}}_0,\widetilde{\mathcal{X}}_0\oplus\widetilde{\mathcal{X}}_0)$
(resp., $\left[\begin{array}{c}
    0 \\
    U
  \end{array}\right]\in
L(\widetilde{\mathcal{X}}_0,\widetilde{\mathcal{X}}_0\oplus\widetilde{\mathcal{X}}_0)$)
on the main diagonal. We will show that the pair
$\mathbf{V}=(V_1,V_2)$ is a uniform isometric $\rho$-dilation of
the pair $\mathbf{A}=(A_1,A_2)$ in the sense of Popescu. First of
all, observe that
\[ V_1^*V_1=I_{\widetilde{\mathcal{X}}}=V_2^*V_2,\quad
V_1^*V_2=V_2^*V_1=0.\] Next, the following relations hold:
\begin{eqnarray*}
\forall k\in\mathbb{N}, \ A_1^k= \left[\begin{array}{cc}
  B^k & 0 \\
  0 & 0
\end{array}\right]=\left[\begin{array}{cc}
  \rho P_{\mathcal{X}_0}U^k|\mathcal{X}_0 & 0 \\
  0 & 0
\end{array}\right] &=& \rho P_\mathcal{X}V_1^k|\mathcal{X};\\
\forall k,n\in\mathbb{N}, \ \forall i_1,\ldots,i_n\in\{
1,2\},\ A_1^kA_2A_{i_1}\cdots A_{i_n} &=& 0\\
&=& \rho P_\mathcal{X}V_1^kV_2V_{i_1}\cdots
V_{i_n}|\mathcal{X}
\end{eqnarray*}
(since $A_1^kA_2=0,\
P_{\widetilde{\mathcal{X}}_0\oplus\widetilde{\mathcal{X}}_0\oplus\{
0\}\oplus\{ 0\}\oplus\ldots}V_1^kV_2=0$);
\begin{eqnarray*}
 A_2= \left[\begin{array}{cc}
  0 & 0 \\
  B & 0
\end{array}\right]=\left[\begin{array}{cc}
  0 & 0 \\
  \rho P_{\mathcal{X}_0}U|\mathcal{X}_0 & 0
\end{array}\right] &=& \rho P_\mathcal{X}V_2|\mathcal{X};\\
\forall k,n\in\mathbb{N}, \ \forall i_1,\ldots,i_n\in\{
1,2\},\ A_2^{k+1}A_{i_1}\cdots A_{i_n} &=& 0\\
&=& \rho P_\mathcal{X}V_2^{k+1}V_{i_1}\cdots
V_{i_n}|\mathcal{X}
\end{eqnarray*}
(since $A_2^2=0,\
P_{\widetilde{\mathcal{X}}_0\oplus\widetilde{\mathcal{X}}_0\oplus\{
0\}\oplus\{ 0\}\oplus\ldots}V_2^2=0$);
\begin{eqnarray*}
\forall k\in\mathbb{N}, \ A_2A_1^k = \left[\begin{array}{cc}
  0 & 0 \\
  B^{k+1} & 0
\end{array}\right] &=& \left[\begin{array}{cc}
  0 & 0 \\
  \rho P_{\mathcal{X}_0}U^{k+1}|\mathcal{X}_0 & 0
\end{array}\right] = \rho P_\mathcal{X}V_2V_1^k|\mathcal{X};\\
\forall k\in\mathbb{N},\ A_2A_1^kA_2=0 &=&
P_\mathcal{X}V_2V_1^kV_2|\mathcal{X};\\
\forall k,n\in\mathbb{N},\ \forall i_1,\ldots,i_n\in\{
1,2\}, &\ &  A_2A_1^kA_2A_{i_1}\cdots A_{i_n}=0\\
&=& \rho P_\mathcal{X}V_2V_1^kV_2V_{i_1}\cdots
V_{i_n}|\mathcal{X}
\end{eqnarray*}
(since $A_1^kA_2=0,\
P_{\widetilde{\mathcal{X}}_0\oplus\widetilde{\mathcal{X}}_0\oplus\{
0\}\oplus\{ 0\}\oplus\ldots}V_2V_1^kV_2=0$). Finally, we get
\[ \forall n\in\mathbb{N},\ \forall i_1,\ldots,i_n\in\{
1,2\}, \   A_{i_1}\cdots A_{i_n}=\rho
P_\mathcal{X}V_{i_1}\cdots V_{i_n}|\mathcal{X}.\]
Thus, $\mathbf{V}$ is a uniform isometric $\rho$-dilation of
$\mathbf{A}$ in the sense of Popescu. However, for any
$\zeta\in\mathbb{T}^N$,
\[ \|\zeta\mathbf{A}\|=\left\|\left[\begin{array}{cc}
  \zeta_1B & 0 \\
  \zeta_2B & 0
\end{array}\right]\right\|=\sqrt{2}\| B\|=\sqrt{2}\rho>\rho.\]
Therefore, $\zeta\mathbf{A}\notin C_\rho$ for all
$\zeta\in\mathbb{T}^N$. We obtain $\mathbf{A}\in
C_{\rho,2}^P\backslash C_{\rho,2}^u$ (moreover, $\mathbf{A}\notin
C_{\rho,2}$).

For the case $N>2$ (and any $\rho>0$) an analogous example of
$\widetilde{\mathbf{A}}\in C_{\rho,N}^P\backslash C_{\rho,N}^u$ is
easily obtained from the previous one, by setting zeros for the
rest of operators in the $N$-tuple, i.e.,
$\widetilde{\mathbf{A}}:=(A_1,A_2,0,\ldots,0)$, where $A_1$ and
$A_2$ are defined in (\ref{eq:a's}). In this case the construction
of a uniform isometric $\rho$-dilation of $\mathbf{A}$ in the
sense of Popescu should be slightly changed (we leave this to a
reader as an easy exercise).
\end{proof}
\begin{rem}
The pair
$\mathbf{A}^{(\varepsilon)}=(A_1^{(\varepsilon)},A_2^{(\varepsilon)})$
constructed in Theorem~\ref{thm:non-similar} doesn't belong to the
class $C_{\rho,2}^u$ for any $\varepsilon>0$ and $\rho>1$. Indeed,
we have shown in Theorem~\ref{thm:non-similar} that
$\mathbf{A}^{(\varepsilon)}$ is not simultaneously similar to any
$\mathbf{T}=(T_1,T_2)\in C_{1,2}$, not speaking of $\mathbf{T}\in
C_{1,2}^u$. Thus, by Theorem~\ref{thm:pop},
$\mathbf{A}^{(\varepsilon)}\notin C_{\rho,2}^u$. This can be shown
also by the following estimate: if $\mathbf{A}^{(\varepsilon)}\in
C_{\rho,2}^u$ for some $\varepsilon>0$ and $\rho>1$, then there
exists a uniform unitary $\rho$-dilation
$\mathbf{U}^{(\varepsilon)}=(U_1^{(\varepsilon)},U_2^{(\varepsilon)})$
of
$\mathbf{A}^{(\varepsilon)}=(A_1^{(\varepsilon)},A_2^{(\varepsilon)})$,
and for any $n\in\mathbb{N}$,
\begin{eqnarray*}
\|
[(A_1^{(\varepsilon)}+A_2^{(\varepsilon)})(A_1^{(\varepsilon)}-A_2^{(\varepsilon)})]^n\|
&=& \|\rho
P_\mathcal{X}[(U_1^{(\varepsilon)}+U_2^{(\varepsilon)})(U_1^{(\varepsilon)}-U_2^{(\varepsilon)})]^n|\mathcal{X}\|
\\
&\leq& \rho\|
[(U_1^{(\varepsilon)}+U_2^{(\varepsilon)})(U_1^{(\varepsilon)}-U_2^{(\varepsilon)})]^n\|=\rho<\infty.
\end{eqnarray*}
This contradicts to (\ref{eq:lim}). Thus, for each $\rho>1$ we
obtain for $\varepsilon>0$ small enough,
$\mathbf{A}^{(\varepsilon)}=(A_1^{(\varepsilon)},A_2^{(\varepsilon)})\in
C_{\rho,2}\backslash C_{\rho,2}^u$, as well as
$\widetilde{\mathbf{A}}:=(A_1,A_2,0,\ldots,0)\in
C_{\rho,N}\backslash C_{\rho,N}^u$.
\end{rem}

\subsection*{Acknowledgements}
I am grateful for the hospitality of the Universities of Leeds and
Newcastle upon Tyne where a part of this work was carried out
during my visits under the International Short Visit Scheme of the
LMS (grant no. 5620). I wish to thank also Dr. Michael Dritschel
from the University of Newcastle upon Tyne for useful discussions.

----------------------------

\begin{thebibliography}{99}% Replace 9 by 99 if 10 or more references
%
% Please note the use of "\and" between author names below
%
\bibitem{Ag}
 {J. Agler},
 `On the representation of certain holomorphic functions defined on a
 polydisc', in Topics in Operator Theory: Ernst D.~Hellinger Memorial
  Volume (L.~de~Branges, I.~Gohberg, and J.~Rovnyak, eds.),
 {\em  Oper.\ Theory\ Adv.\ Appl.\ }48 (1990) 47--66 (Birkh\"auser-Verlag,
 Basel).
 %
\bibitem{AgMc}
{J.~Agler \and J.~E.~McCarthy}, `Nevanlinna--Pick
interpolation on the bidisk', {\em  J.\ Reine\ Angew.\ Math.\ }
506 (1999) 191--204.
%
\bibitem{An}
{T.~Ando}, `On a pair of commutative contractions', {\em
Acta\ Sci.\ Math.\ (Szeged)\ } 24 (1963) 88--90.
%
\bibitem{AnN}
{T.~Ando \and K.~Nishio}, `Convexity properties of
operator radii associated with unitary $\rho$-dilations', {\em
Michigan\ Math.\ J.\ } 20 (1973) 303--307.
%
\bibitem{BaCa}
{C.~Badea \and G.~Cassier}, `Constrained von Neumann
inequalities', {\em  Adv.\ Math.\ } 166\ no.~2 (2002) 260--297.
%
\bibitem{BLTT}
{J.~A.~Ball, W.~S.~Li, D.~Timotin \and T.~T.~Trent}, `A
commutant lifting theorem on the polydisc', {\em Indiana\ Univ.\
Math.\ J.\ } 48\ no.~2 (1999) 653--675.
%
\bibitem{BSV}
{J.~A.~Ball, C.~Sadosky \and V.~Vinnikov}, `Conservative
input-state-output systems with evolution on a multidimensional
integer lattice', preprint.
%
\bibitem{BT}
{J.~A.~Ball \and T.~T.~Trent}, `Unitary colligations,
reproducing kernel Hilbert spaces, and Nevanlinna-Pick
interpolation in several variables', {\em J.\ Funct.\ Anal.\ }
157\ no.~1 (1998) 1--61.
%
\bibitem{Be}
{C.~A.~Berger}, `A strange dilation theorem', {\em
Notices\ Amer.\ Math.\ Soc.\ } 12 (1965) 590.
%
\bibitem{CaF}
{G.~Cassier \and T.~Fack}, `Contractions in von Neumann
algebras', {\em  J.\ Funct.\ Anal.\ } 135\ no.~2 (1996) 297--338.
%
\bibitem{Da}
{C.~Davis}, `The shell of a Hilbert-space operator', {\em
Acta\ Sci.\ Math.\ (Szeged)\ }  29 (1968) 69--86.
%
\bibitem{DMcCW}
{M.~A.~Dritschel, S.~McCullough \and H.~J.~Woerdeman},
`Model theory for $\rho$-contractions, $\rho\le2$', {\em  J.\
Operator\ Theory\ } 41\ no.~2 (1999) 321--350.
%
\bibitem{Du}
{E.~Durszt}, `On unitary $\rho$-dilations of operators',
{\em Acta\ Sci.\ Math.\ (Szeged)\ } 27 (1966) 247--250.
%
\bibitem{FH}
{C.~K.~Fong \and J.~A.~R.~Holbrook}, `Unitarily invariant
operator norms', {\em  Canad.\ J.\ Math.\ } 35\ no.~2 (1983)
274--299.
%
\bibitem{H1}
{J.~A.~R.~Holbrook}, `On the power-bounded operators of
Sz.-Nagy and Foia\c{s}', {\em Acta\ Sci.\ Math.\ (Szeged)\ } 29
(1968) 299--310.
%
\bibitem{H2}
{J.~A.~R.~Holbrook}, `Inequalities governing the operator
radii associated with unitary $\rho $-dilations', {\em  Michigan\
Math.\ J.\ } 18 (1971) 149--159.
%
\bibitem{K2}
{D.~S.~Kalyuzhniy}, `Multiparametric dissipative linear
stationary dynamical scattering systems: discrete case', {\em J.\
Operator\ Theory} 43\ no.~2 (2000) 427--460.
%
\bibitem{K3}
{D.~S.~Kalyuzhniy}, `Multiparametric dissipative linear
stationary dynamical scattering systems: discrete case. II.
Existence of conservative dilations', {\em Integral\ Equations\
Operator\ Theory} 36\ no.~1 (2000) 107--120.
%
\bibitem{K4}
{D.~S.~Kalyuzhniy}, `On the notions of dilation,
  controllability, observability, and minimality in the theory of
  dissipative scattering linear nD systems', in Proceedings CD of the
   Fourteenth International Symposium of Mathematical Theory of
Networks and Systems (MTNS), June 19--23, 2000, Perpignan, France
(A.~El~Jai and M.~Fliess, Eds.), or \verb
http://www.univ-perp.fr/mtns2000/articles/SI13_3.pdf/.
%
\bibitem{K6}
{D.~S.~Kalyuzhniy-Verbovetzky}, `Cascade connections of
linear systems and factorizations of holomorphic operator
functions around a multiple zero in  several variables', {\em
Math.\ Rep.\ (Bucur.)} 3(53)\ no.~4 (2001) 323--332.
%
\bibitem{K7}
{D.~S.~Kalyuzhniy-Verbovetzky}, `On $J$-conservative
scattering system realizations in several variables', {\em
Integral\ Equations\ Operator\ Theory } 43\ no.~4 (2002) 450--465.
%
\bibitem{K1}
{D.~S.~Kalyuzhny\u\i}, `The von Neumann inequality for
linear matrix functions of several variables', {\em Mat.\ Zametki}
64\ no.~2 (1998) 218--223 (Russian). English transl. in {\em
Math.\ Notes } 64\ no.~1--2 (1998) 186--189 (1999).
%
\bibitem{K5}
{D.~S.~Kalyuzhny\u\i-Verbovetski\u\i}, `Cascade
connections of multiparameter linear systems and the conservative
realization of a decomposable inner operator function on the
bidisk', {\em Mat.\ Stud. } 15\ no.~1 (2001) 65--76 (Russian).
%
\bibitem{NO}
{T.~Nakazi \and K.~Okubo}, `$\rho$-contraction and
$2\times 2$ matrix', {\em  Linear\ Algebra\ Appl.\ } 283\ no.~1-3
(1998) 165--169.
%
\bibitem{vN}
{J.~von Neumann}, `Eine Spektraltheorie f\"{u}r
allgemeine Operatoren eines unit\"{a}ren Raumes', {\em Math.\
Nachr. } 4 (1951) 258--281 (German).
%
\bibitem{OAn}
{K.~Okubo \and T.~Ando}, `Operator radii of commuting
products', {\em  Proc.\ Amer.\ Math.\ Soc.\ } 56\ no.~1 (1976)
203--210.
%
\bibitem{OS}
{K.~Okubo \and I.~Spitkovsky}, `On the characterization
of $2\times 2$ $\rho$-contraction matrices', {\em  Linear\
Algebra\ Appl.\ } 325\ no.~1-3 (2001) 177--189.
%
\bibitem{Pa}
{V.~I.~Paulsen}, `Every completely polynomially bounded
operator is similar to a contraction', {\em J.\ Funct.\ Anal.\ }
55\ no.~1 (1984) 1--17.
%
\bibitem{Po1}
{G.~Popescu}, `Isometric dilations for infinite sequences
of noncommuting operators', {\em  Trans.\ Amer.\ Math.\ Soc.\ }
316\ no.~2 (1989) 523--536.
%
\bibitem{Po2}
{G.~Popescu}, `Positive-definite functions on free
semigroups', {\em  Canad.\ J.\ Math.\ } 48\ no.~4 (1996) 887--896.
%
\bibitem{SzN}
{B.~Sz.-Nagy}, `Sur les contractions de l'espace de
Hilbert', {\em  Acta\ Sci.\ Math.\ (Szeged)\ } 15 (1953) 87--92
(French).
%
\bibitem{SzNF1}
{B.~Sz.-Nagy \and C. Foia\c{s}}, `On certain classes of
power-bounded operators in Hilbert space', {\em  Acta\ Sci.\
Math.\ (Szeged)\ } 27 (1966) 17--25.
%
\bibitem{SzNF2}
{B.~Sz.-Nagy \and C. Foia\c{s}}, `Similitude des
op\'{e}rateurs de class ${\mathcal C}\sb{\rho }$ \`{a} des
contractions', {\em C.\ R.\ Acad.\ Sci.\ Paris\ S\'{e}r.\ A-B\ }
264 (1967) A1063--A1065 (French).
%
\bibitem{SzNF}
 {B.~Sz.-Nagy \and C. Foia\c{s}},
 {\em Harmonic analysis of operators on Hilbert space}
 (North-Holland, Amsterdam--London, 1970).
%
\bibitem{W}
{J.~P.~Williams}, `Schwarz norms for operators', {\em
Pacific\ J.\ Math.\ } 24 (1968) 181--188.
%

\end{thebibliography}
\end{document}